\documentstyle[twoside,11pt]{article}
\begin{document}

\title{\bf One-dimensional Compressible Navier-Stokes Equations with Temperature Dependent Transport Coefficients and Large Data} \vskip 0.2cm
\author{
{\bf Hongxia Liu}\\
Department of Mathematics, Jinan University\\
Guangzhou 510632, China\\[3mm]
{\bf Tong Yang}\\
Department of Mathematics, City University of Hong Kong\\
Tat Chee Avenue, Kowloon, Hong Kong, China\\[3mm]
{\bf Huijiang Zhao}\\
School of Mathematics and Statistics\\
Wuhan University, Wuhan 430072,  China\\[3mm]
{\bf Qingyang Zou}\\
School of Mathematics and Statistics\\
Wuhan University, Wuhan 430072,  China}
\date{ }

\maketitle

\vskip 0.2cm
\arraycolsep1.5pt
\newtheorem{Lemma}{Lemma}[section]
\newtheorem{Theorem}{Theorem}[section]
\newtheorem{Definition}{Definition}[section]
\newtheorem{Proposition}{Proposition}[section]
\newtheorem{Remark}{Remark}[section]
\newtheorem{Corollary}{Corollary}[section]

\begin{abstract}
This paper is concerned with the global smooth non-vacuum solutions with large data to the Cauchy problem of the one-dimensional compressible Navier-Stokes equations with degenerate temperature dependent transport coefficients which satisfy conditions from the consideration in kinetic theory. A Nishida-Smoller type result is obtained.\\

\noindent{\bf A.M.S. Mathematics Subject Classification:} 35Q35, 35D35, 74D10, 76D05, 76N10
\end{abstract}

\tableofcontents

\section{Introduction and main result}
\setcounter{equation}{0}

The motion of one-dimensional compressible flow of a  viscous ideal fluid can be described by the  Navier-Stokes system:
\begin{eqnarray}\label{NS-Euler}
\left\{
\begin{array}{rl}
\rho_\tau+(\rho u)_y=&0,\\[2mm]
(\rho u)_\tau+\left(\rho u^2+p(\rho,\theta)\right)_y=&\left(\mu u_y\right)_y,\\[2mm]
(\rho{\mathcal{E}})_\tau+\left(\rho u{\mathcal{E}}+up(\rho,\theta)\right)_y=&\left(\kappa \theta_y+\mu uu_y\right)_y.
\end{array}
\right.
\end{eqnarray}
Here $y$ and $\tau$ represent  the space variable and the time variable
respectively, and the primary dependent variables are fluid density $\rho$, fluid velocity $u$, and temperature $\theta$. The specific total energy ${\mathcal{E}}=e+\frac 12 u^2$ with $e$ being the specific internal energy. The pressure $p$, the internal energy
$e$, and the transport coefficients $\mu>0$ (viscosity) and $\kappa>0$ (heat conductivity) are  functions of $\rho$ and $\theta$. The thermodynamic variables $\rho, p, e, s$, and $\theta$ are related by the Gibbs equation $de =¦Èds-pd\rho^{-1}$, where $s$ is the specific entropy.

Motivated by the study in the kinetic theory of gases,
we are interested in  constructing  global smooth non-vacuum solutions to the Cauchy problem of the system (\ref{NS-Euler}) with large initial data for the case when the transport coefficients $\mu>0$ and $\kappa>0$ are functions of temperature.

More precisely, recall that the Boltzmann equation with  slab symmetry takes the form
\begin{equation}\label{Boltzmann}
f_\tau+\xi_1f_y=\frac{1}{\varepsilon}Q(f,f),
\end{equation}
where the unknown function $f(\tau,y,\xi)\geq 0$ stands for the distribution density of particles with position $y\in {\bf R}$ and velocity $\xi=(\xi_1,\xi_2,\xi_3)\in{\bf R}^3$ at time $\tau\geq 0$, $\varepsilon>0$  is the Knudsen number proportional to the mean free path and it measures the adiabaticity of the gas, and $Q(f,f)$ is the Boltzmann collision operator defined by
$$
Q(f,g)=\frac 12\int_{{\bf R}^3\times {\bf S}^2}q(|\xi-\xi_*|,\theta)\Big(f'g'_*+f'_*g'-fg_*-f_*g\Big)d\xi_*d\omega.
$$
Here $q(|\xi-\xi_*|,\theta)\geq 0$ is the cross section that
 is determined by the interaction potential of two colliding particles.
Here,  $f=f(\tau,y,\xi),\ f_*=f(\tau,y,\xi_*),\
f'=f(\tau,y,\xi'),\ f'_*=f(\tau,y,\xi'_*)$,
$\cos\theta=(\xi-\xi_*)\cdot\omega/|\xi-\xi_*|,\ \omega\in{\bf S}^2$,
and
$$
\xi'=\xi-((\xi-\xi_*)\cdot\omega)\omega,\quad \xi'_*=\xi_*+((\xi-\xi_*)\cdot\omega)\omega,
$$
is the relation between the velocities $\xi',\ \xi'_*$ after and the velocities $\xi,\ \xi_*$ before the collision by the conservation of momentum and energy. For details see \cite{Cercignani-Illner-Pulvirenti}.

It is well-known that, by employing the celebrated Chapman-Enskog expansion, the compressible Navier-Stokes equations (\ref{NS-Euler}) is the first order
 approximation of the Boltzmann equation (\ref{Boltzmann}) in term of
$\epsilon$ and  the viscosity $\mu$ and heat conductivity $\kappa$ are functions of temperature,  cf. \cite{Cercignani-Illner-Pulvirenti}, \cite{Chapman-Cowling}, \cite{Grad}, \cite{Vincenti-Kruger}. In particular, if the inter-molecule potential is propotional to $r^{-\alpha}$ with $r$ being the molecule distance. Then
$$
\mu,\ \ \kappa\propto \theta^{\frac{\alpha+4}{2\alpha}}.
$$
Note that for Maxwellian molecules ($\alpha$=4) the dependence is linear, while for elastic spheres $(\alpha\rightarrow+\infty)$ the dependence is like $\sqrt{\theta}$.

The above dependence has strong influence on the solution behavior and leads to difficulty in analysis for global existence with large data.
In fact, as pointed out in \cite{Jenssen-Karper}, temperature dependence of the viscosity $\mu$  turns out be especially problematic and challenging.  Even though
there are works about the density dependence in $\mu$ and temperature and density dependence in $\kappa$,  cf. \cite{Antontsev-Kazhikhov-Monakhov} \cite{Dafermos}, \cite{Jenssen-Karper}, \cite{Jiang-Racke},  \cite{Kanel},  \cite{Kawashima-Nishida}, \cite{Kawashima-Okada}, \cite{Kawohl}, \cite{Kazhikhov}, \cite{Kazhikhov-Shelukhin}, \cite{Tan-Yang-Zhao-Zou} and the references therein,
 no result has been obtained for the case when $\mu$ depends on temperature. Hence, the result obtained in the paper can be viewed a small progress in this direction.

Throughout this paper, we will assume
\begin{equation}\label{1.4}
\mu=\mu(\theta)>0,\quad \kappa=\kappa(\theta)>0,\quad \forall \theta>0.
\end{equation}

To state the main result, let $x$ be the Lagrangian space variable, $t$ be the time variable, and $v=\frac 1\rho$ denote the specific volume. Then the
system  (\ref{NS-Euler}) becomes
\begin{eqnarray}\label{1.5}
\left\{
\begin{array}{rl}
v_t-u_x=&0,\\[2mm]
u_t+p(v,\theta)_x=&\left(\frac{\mu(\theta)u_x}{v}\right)_x,\\[2mm]
e_t+p(v,\theta)u_x=&\frac{\mu(\theta)u_x^2}{v}
+\left(\frac{\kappa(\theta)\theta_x}{v}\right)_x.
\end{array}
\right.
\end{eqnarray}
Assume that  the initial data satisfy
\begin{equation}\label{1.6}
(v(0,x), u(0,x), \theta(0,x))=(v_0(x), u_0(x), \theta_0(x)),\quad \lim\limits_{x\to\pm\infty}(v_0(x), u_0(x), \theta_0(x))=(1, 0,1).
\end{equation}

Throughout this paper, we will concentrate on the ideal, polytropic gases:
\begin{equation}\label{1.7}
p(v,\theta)=\frac{R\theta}{v}=Av^{-\gamma}\exp\left(\frac{\gamma-1}{R}s\right), \quad
e=C_v\theta=\frac{R\theta}{\gamma-1},
\end{equation}
where the specific gas constants $A$, $R$ and the specific heat at constant volume $C_v$ are positive constants and $\gamma>1$ is the adiabatic constant.
Note that for the model of monatomic gas, $\gamma=5/3$, that does not
satisfy the condition imposed in the following theorem.

We will present a Nishida-Smoller type result for the above problem (For the corresponding Nishida-Smoller type global existence result for one-dimensional ideal polytropic compressible Euler system, please refer to \cite{Nishida-Smoller}). To state the main result, we will choose the velocity $u$,  the specific volume $v$, and the entropy $s$ as unknown functions, and  use $\overline{s}=\frac{R}{\gamma-1}\ln\frac{R}{A}$ to denote the far field of the initial entropy $s_0(x)$, that is,
$$
\lim\limits_{|x|\to+\infty}s_0(x)=\lim\limits_{|x|\to+\infty}\frac{R}{\gamma-1}
\ln{\frac{R\theta_0(x)v_0(x)^{\gamma-1}}{A}=\bar s}.
$$

\begin{Theorem} Suppose that
\begin{itemize}
\item[$\bullet$]  $N_0=\|(v_0(x)-1,u_0(x),s_0(x)-\bar{s})\|_{H^3({\bf R})}$ is bounded by some positive constant independent of $\gamma-1$ and there are $(\gamma-1)-$independnet positive constants $0<\underline{V}_0<1,\ \overline{V}_0>1,\ 0<\underline{\Theta}_0<1,\ \overline{\Theta}_0>1$ such that
$$
\underline{V}_0\leq v_0(x)\leq \overline{V}_0,~~\underline{\Theta}_0\leq \theta_0(x)\leq \overline{\Theta}_0, \quad \forall x\in{\bf R};
$$

\item[$\bullet$] $\mu(\theta)$ and $\kappa(\theta)$ are smooth for $\theta>0$ and satisfy (\ref{1.4}) for $\theta>0$;
\item[$\bullet$] There exists a non-negative continuous function $C(N_0)$ satisfying $C(x)>0$ for $x>0$ and $C(0)=0$ such that
\begin{equation}\label{1.8}
(\gamma-1)C(N_0)\leq 1.
\end{equation}
\end{itemize}
Then the Cauchy problem (\ref{1.5}), (\ref{1.6}) admits a unique global solution $(v(t,x), u(t,x),\theta(t,x))$ satisfying
\begin{equation}\label{1.9}
\underline{V}\leq v(t,x)\leq \overline{V},\quad \underline{\Theta}_0\leq \theta(t,x)\leq \overline{\Theta}_0,
\end{equation}
 and
\begin{equation}\label{1.10}
\lim\limits_{t\to+\infty}\sup\limits_{x\in{\bf R}}\Big|(v(t,x)-1,u(t,x),\theta(t,x)-1)\Big|=0.
\end{equation}
Here $\underline{V}$ and $\overline{V}$ are some positive constants  depending only on $N_0$, $\underline{V}_0$, $\overline{V}_0$, $\underline{\Theta}_0$, and $\overline{\Theta}_0$.
\end{Theorem}

\begin{Remark} Several remarks concerning Theorem 1.1 are given:
\begin{itemize}
\item Even for the case when the viscosity coefficient $\mu$ and the heat conductivity coefficient $\kappa$ are functions of both $v$ and $\theta$, similar result still holds if $\mu(v,\theta)$ and $\kappa(v,\theta)$ satisfy
    \begin{itemize}
    \item[$\bullet$] $\mu(v,\theta)>0, \kappa(v,\theta)>0$ hold for $v>0,\ \theta>0$,
    \item [$\bullet$] There exist constants $a<0, \ b>-\frac 12$ such that
    $$
    \mu(v,\theta)\sim\left\{
    \begin{array}{rl}
    v^a,\quad & v\to 0_+,\\[2mm]
    v^b,\quad & v\to+\infty,
    \end{array}
    \right.
    $$
     for any $\theta\in \left[\underline{\Theta}_0,\overline{\Theta}_0\right]$.
    \end{itemize}

\item From the proof of Theorem 1.1, it is easy to see that even when $\|(v_0-1,u_0,\theta_0-1)\|_{H^3({\bf R})}$, $\underline{V}_0$, $\overline{V}_0$, $\underline{\Theta}_0$, and $\overline{\Theta}_0$ depend on $\gamma-1$, with
    $\lim\limits_{\gamma\to 1_+}\|(v_0-1,u_0,\theta_0-1)\|_{H^3({\bf R})}=+\infty$, $\lim\limits_{\gamma\to 1_+}\left(\underline{V}_0,\underline{\Theta}_0\right)=(0,0)$,
    $\lim\limits_{\gamma\to 1_+}\left(\overline{V}_0,\overline{\Theta}_0\right)=(+\infty, +\infty),$ a similar result as Theorem 1.1 also holds provided that the above limits satisfy certain growth conditions when $\gamma\to 1_+$.
\end{itemize}
\end{Remark}

Now we sketch the main ideas used in the proof of Theorem 1.1. As pointed out in \cite{Antontsev-Kazhikhov-Monakhov, Dafermos, Kawohl, Tan-Yang-Zhao-Zou}, the key point to the global solvability of the Cauchy problem (\ref{1.5}), (\ref{1.6}) with large data is to obtain the positive lower and upper bounds for the specific volume $v(t,x)$ and the absolute temperature $\theta(t,x)$. For the Cauchy problem (\ref{1.5}), (\ref{1.6}), to  our knowledge, according to the dependence of
the viscosity coefficient $\mu$ and the heat conductivity coefficient $\kappa$  on $v$ and $\theta$, there are two existing effective approaches:
\begin{itemize}
\item[(i).] The first approach is developed by A. V. Kazhikhov and V. V. Shelukhin in \cite{Kazhikhov-Shelukhin} for the case when the viscosity coefficient $\mu$ is a positive constant, i.e., $\mu=\mu_0>0$ for some constant $\mu_0$. In fact, only the case when both $\mu$ and $\kappa$ are positive constants is discussed in \cite{Kazhikhov-Shelukhin}. However, the argument developed in \cite{Kazhikhov-Shelukhin} can be applied to  the case when $\mu$ is a positive constant and $\kappa$
 depends on both $v$ and $\theta$. The main idea in \cite{Kazhikhov-Shelukhin} is to deduce an explicit formula for $v(t,x)$, cf. \cite{Kazhikhov-Shelukhin}, as follows:

 For each $i\in{\bf Z}$ and $x\in[i,i+1]$
\begin{equation}\label{1.11}
v(t,x)=\frac{1+\frac {R}{\mu_0} {\displaystyle\int^t_0}\theta(\tau,x)B_i(\tau,x)Y_i(\tau)d\tau}{B_i(t,x)Y_i(t)}.
\end{equation}
Here
\begin{eqnarray*}
B_i(t,x)&=&\frac{v_0(a_i(t))}{v_0(x)v(t,a_i(t))}\exp\left(\frac {1}{\mu_0} \int^{a_i(t)}_x(u(t,y)-u_0(y))dy\right),\\
Y_i(t)&=&\exp\left(\frac {R}{\mu_0}\int^t_0\left(\frac\theta v\right)(\tau, a_i(t))d\tau\right),
\end{eqnarray*}
and for each integer $i\in{\bf Z}$, $a_i(t)\in[i,i+1]$ satisfies
$$
A_1\leq v(t,a_i(t))\leq A_2,
$$
with $A_i>0$ being the two positive roots of $x-\ln x+1=C$ for some sufficiently large positive constant $C$ depending only on the initial data.

With the expression (\ref{1.11}),  one can firstly
deduce  a positive lower bound for $v(t,x)$. And then by employing the standard maximum principle for the following parabolic equation
\begin{eqnarray}\label{1.12}
C_v\left(\frac 1\theta\right)_t&=&
-\frac{\mu u_x^2}{\theta^2v}+\frac{Ru_x}{v\theta}-
   \frac{2\theta\kappa}{v}\left[\left(\frac1\theta\right)_x\right]^2
       +\left[\left(\frac{\kappa}v\right)\left(\frac1\theta\right)_x\right]_x \nonumber\\
   &=&\left[\left(\frac{\kappa}v\right)\left(\frac1\theta\right)_x\right]_x
      -\left\{\frac{2\theta\kappa}{v}\left[\left(\frac1\theta\right)_x\right]^2
      +\frac{\mu}{v\theta^2}
        \left(u_x-\frac{R\theta}{2\mu} \right)^2\right\} \\
     &&+\frac {R^2}{4v\mu},\nonumber
  \end{eqnarray}
one can  deduce a positive lower bound for $\theta(t,x)$.

With the lower bounds on both $v$ and $\theta$, the argument used in \cite{Kazhikhov-Shelukhin} leads to the upper bound estimate on $v$ provided that $\kappa$ satisfies
$$
\min\limits_{v\geq V_1>0, \theta\geq \Theta_1>0}\kappa(v,\theta)\geq C\left(V_1,\Theta_1\right)>0
$$
for some positive constant $C\left(V_1,\Theta_1\right)>0$. And then
the upper bound on $\theta$ follows and  the global existence of solution is proved.

\item[(ii).] The second approach was introduced in \cite{Tan-Yang-Zhao-Zou} to treat the case when $\mu$ and $\kappa$ are degenerate functions of $v$ and/or $\theta$, say for example $\mu=v^{-a},$ $\kappa=\theta^b$ for some positive constants $a>0,\ b>0$. In such a case, the argument used in \cite{Kazhikhov-Shelukhin} can not be used. And the main idea in \cite{Tan-Yang-Zhao-Zou} is to firstly derive
  the lower bound for $\theta$ in term of the lower bound of $v$
\begin{equation}\label{1.13}
\left\|\frac 1\theta\right\|_{L^\infty([0,T]\times{\bf R})}\leq O(1)\left\{1+
\left\|\frac 1 v\right\|^{1-a}_{L^\infty([0,T]\times{\bf R})}\right\},\quad 0\leq a<1,
\end{equation}
by applying the maximum principle to (\ref{1.12}), and then to deduce the following lower and upper bounds for $v$ in terms of $\theta$
\begin{equation}\label{1.14}
\left\|\frac1{v}\right\|_{L^\infty([0,T]\times{\bf R})}\leq O(1)\left(1+\left\|\theta^{1-b}\right\|_{L^\infty([0,T]\times{\bf R})}^{\frac1{3a-1}}\right),\quad \frac 13<a<\frac 12,
\end{equation}
and
\begin{equation}\label{1.15}
\|v\|_{L^\infty([0,T]\times{\bf R})}\leq O(1)\left(1+\left\|\theta^{1-b}\right\|_{L^\infty{([0,T]\times{\bf R})}}^{\frac{2a}{(3a-1)(1-2a)}}\right),\quad \frac 13<a<\frac 12,
\end{equation}
by applying Y. Kanel's argument \cite{Kanel}.

From the above estimates, we have
\begin{eqnarray}\label{1.16}
\|\theta\|_{L^\infty([0,T]\times{\bf R})}&\leq& O(1)\left\{1+\int_0^t\left(\left\|\frac{u_x^2}{v^{1+a}}\right\|_{L^\infty([0,T]\times{\bf R})}\right.\right.\\
&&\left.\left.+
 \left\|\frac{u_x^2}{v^{2}}\right\|_{L^\infty([0,T]\times{\bf R})}+\|\theta\|^2_{L^\infty([0,T]\times{\bf R})}\right)d\tau\right\}.\nonumber
\end{eqnarray}
With this relation, one can deduce the desired lower and upper bounds on $v$ and $\theta$ if $a$ and $b$ satisfy certain conditions.
\end{itemize}

We note that in the above analysis one can obtain from (\ref{1.5})$_1$ and (\ref{1.5})$_2$ that
\begin{equation}\label{1.17}
\left(\frac{\mu(v)v_x}{v}\right)_t=\left(\frac{\mu(v)v_t}{v}\right)_x
=\left(\frac{\mu(v)u_x}{v}\right)_x=u_t+p(v,\theta)_x.
\end{equation}
From which and the fact that the gas under consideration is ideal polytropic, one can then deduce an estimate on $\left\|\frac{\mu(v)v_x}{v}\right\|$.

When the viscosity coefficient $\mu$ depends on $\theta$, the situation is different because the identity corresponding to (\ref{1.17}) now becomes
\begin{equation}\label{1.18}
\left(\frac{\mu(v,\theta)v_x}{v}\right)_t=u_t+p(v,\theta)_x+
\frac{\mu_\theta(v,\theta)}v\left(\theta_tv_x-u_x\theta_x\right),\quad \mu_\theta(v,\theta)=\frac{\partial \mu(v,\theta)}{\partial\theta}.
\end{equation}
From  (\ref{1.5})$_3$,  the last term on the right hand side of (\ref{1.18}) is  highly nonlinear so that to control the possible growth of $(v(t,x), u(t,x), \theta(t,x))$  is very difficult for  large initial perturbation.

The key observation  in this paper is  that the constitutive relations (\ref{1.7}) gives
$$
\theta(t,x)=\frac AR 
v(t,x)^{1-\gamma}\exp
\left(\frac{\gamma-1}{R}(s(t,x)-\bar{s})\right)\exp\left(\frac{\gamma-1}{R}\bar{s}\right)
$$
which implies that the absolute temperature $\theta(t,x)$ can be sufficiently close to $1$ when $\gamma-1>0$ is sufficiently small, and $v(t,x)$ is bounded from both below and above by some positive constants independent of $\gamma-1$ together with
 $\lim\limits_{\gamma\to 1_+}(\gamma-1)\|s(t,x)-\bar{s}\|_{L^\infty([0,T]\times{\bf R})}=0$.
Thus, under the {\it a priori assumption} on the absolute temperature $\theta(t,x)$
\begin{equation}\label{1.19}
\frac 12\underline{\Theta}_0\leq \theta(t,x)\leq2\overline{\Theta}_0,
\end{equation}
 by some delicate energy type estimates and  using the argument initiated in \cite{Kanel}, we can deduce an uniform in time positive lower and upper bound on $v(t,x)$ and some uniform energy estimates on $\left\|\left(v-1,u,\frac{\theta-1}{\sqrt{\gamma-1}}\right)(t)\right\|_{H^3({\bf R})}$ in terms of $\left\|\left(v_0-1,u_0,\frac{\theta_0-1}{\sqrt{\gamma-1}}\right)\right\|_{H^3({\bf R})}$, $\inf\limits_{x\in{\bf R}}v_0(x)$, $\sup\limits_{x\in{\bf R}}v_0(x)$, $\inf\limits_{x\in{\bf R}}\theta_0(x)$, and $\sup\limits_{x\in{\bf R}}\theta_0(x)$. These estimates are sufficient to show global existence when $\gamma-1$ is small.

The rest of this paper is organized as follows. In Section 2, we will give some identities
for later use. The energy estimates are given  in Section 3. And
the proof of the main result will be given in Section 4.

\vskip 2mm
\noindent{\bf Notations:}  O(1) or $C_i (i\in {\bf N})$ stands for a generic positive constant which is independent of $t$, $x$, and $\gamma-1$, while $C(\cdot,\cdots,\cdot)$  is used to denote some positive constant depending only on the arguments listed in the parenthesis. Note that all these constants may vary from line to line.
$\|\cdot\|_s$ represents the norm in $H^s({\bf R})$ with $\|\cdot\|=\|\cdot\|_0$ and for $1\leq p\leq +\infty$, $L^p({\bf R})$ denotes the standard Lebesgue space.

Finally, to simplify the presentation, we can assume without loss of generality that the gas constants $A=R=1$ and consequently $\bar{s}=0$.

\section{Preliminaries}
\setcounter{equation}{0}

This section is devoted to listing some identities which will be used in the following sections.

Firstly, notice that
\[
\left(\frac{\mu(\theta)v_x}{v}\right)_t=\left(\frac{\mu(\theta)u_x}{v}\right)_x+
\frac{\mu'(\theta)}v\left(\theta_tv_x-u_x\theta_x\right),
\]
we have from (\ref{1.5})$_2$ that
\begin{equation}\label{2.1}
\left(\frac{\mu(\theta)v_x}{v}\right)_t=u_t+\left(\frac{\theta}{v}\right)_x+
\frac{\mu'(\theta)}v\left(\theta_tv_x-u_x\theta_x\right).
\end{equation}
Here recall that all the gas constants $A$ and $R$ have been normalized to be $1$.

On the other hand, we can get from (\ref{1.5})$_3$ that
\begin{eqnarray}\label{2.2}
\frac{\theta_t}{\gamma-1}&=&\frac{\mu(\theta)u_x^2}{v}
+\left(\frac{\kappa(\theta)\theta_x}{v}\right)_x
-\frac{\theta u_x}{v}\nonumber\\
&=&\frac{\mu(\theta)u_x^2
 +\kappa'(\theta)\theta_x^2+\kappa(\theta)\theta_{xx}-\theta u_x}{v}
 -\frac{\kappa(\theta)\theta_xv_x}{v^2}\nonumber\\
&\leq& C(v,\theta)\left(u_x^2+\theta_x^2+|u_x|+|\theta_xv_x|+|\theta_{xx}|\right),
\end{eqnarray}
\begin{eqnarray}\label{2.3}
\frac{\theta_{tx}}{\gamma-1}&=&\frac{\kappa(\theta)\theta_{xxx}
+3\kappa'(\theta)\theta_x\theta_{xx}+\kappa''(\theta)\theta^3_x}{v}\nonumber\\
&&+\frac{\mu'(\theta)\theta_xu_x^2+2\mu(\theta)u_xu_{xx}-\theta_xu_x-\theta u_{xx}
}{v}\\
&&-\frac{2\kappa(\theta)\theta_{xx}v_x+2\kappa'(\theta)\theta_x^2v_x
+\kappa(\theta)\theta_xv_{xx}}{v^2}\nonumber\\
&&-\frac{\mu(\theta)u_x^2v_x-\theta u_xv_x}{v^2}
+\frac{2\kappa(\theta)\theta_xv_x^2}{v^3}\nonumber\\
&\leq&C(v,\theta)\Big(|\theta_{xxx}|+\left|(v_x,\theta_x)\right|\left|\theta_{xx}\right|+
\left(1+|u_x|\right)\left|u_{xx}\right|+
|\theta_x|\left|v_{xx}\right|\nonumber\\
&&+\left(1+\left|(u_x,\theta_x)\right|\right)
\left|(v_x,u_x,\theta_x)\right|^2\Big),\nonumber
\end{eqnarray}
\begin{eqnarray}\label{2.4}
\frac{\theta_{txx}}{\gamma-1}
&=&\frac{2\kappa(\theta)\theta_xv_x^3}{v^4}
+\frac{6\kappa(\theta)\theta_{xx}v^2_x+6\kappa'(\theta)\theta_x^2v^2_x}{v^3}\nonumber\\
&&+\frac{6\kappa(\theta)\theta_xv_xv_{xx}+2\mu(\theta)u_x^2v^2_x-2\theta u_xv^2_x}{v^3}\nonumber\\
&&-\frac{2\mu'(\theta)u_x^2v_x\theta_x+4\mu(\theta)u_xv_xu_{xx}+
\mu(\theta)u_x^2v_{xx}+3\kappa''(\theta)\theta_x^3v_x}{v^2}\\
&&-\frac{9\kappa'(\theta)\theta_x\theta_{xx}v_x
+3\kappa'(\theta)\theta_x^2v_{xx}+3\kappa(\theta)\theta_{xxx}v_x+3\kappa(\theta)\theta_{xx}v_{xx}}{v^2}\nonumber\\
&&-\frac{\kappa(\theta)\theta_xv_{xxx}
-2\theta_xv_xu_x-2\theta u_{xx}v_x-\theta u_xv_{xx}}{v^2}\nonumber\\
&&+\frac{\mu''(\theta)u^2_x\theta^2_x+\mu'(\theta)u^2_x\theta_{xx}+
4\mu'(\theta)\theta_xu_xu_{xx}}{v}\nonumber\\
&&+\frac{2\mu(\theta)u_{xx}^2+2\mu(\theta)u_{xxx}u_x+\kappa'''(\theta)\theta^4_x
}{v}\nonumber\\
&&+\frac{6\kappa''(\theta)\theta_x^2\theta_{xx}+3\kappa'(\theta)\theta^2_{xx}+4\kappa'(\theta)\theta_{xxx}\theta_{x}
+\kappa(\theta)\theta_{xxxx}}{v}\nonumber\\
&&-\frac{\theta_{xx}u_x+2\theta_xu_{xx}+\theta u_{xxx}}{v}\nonumber\\
&\leq&C(v,\theta)\Big(\left|\theta_{xxxx}\right|+\left|u_{xxx}\right|+\left|(u_x,\theta_x)\right|
\left|\left(v_{xxx},\theta_{xxx}\right)\right|+\left|\left(v_{xx},u_{xx},\theta_{xx}\right)\right|^2\nonumber\\
&&+\left|\left(v_{xx},u_{xx},\theta_{xx}\right)\right|\left|\left(v_{x},u_{x},\theta_{x}\right)\right|
\left(1+\left|\left(v_{x},u_{x},\theta_{x}\right)\right|\right)
+\left|\left(v_{x},u_{x},\theta_{x}\right)\right|^3
\left(1+\left|\left(v_{x},u_{x},\theta_{x}\right)\right|\right)\Big)\nonumber.
\end{eqnarray}
Next, we give some identities related to the pressure $p$. Recall that $p=\frac\theta v$, we have
\begin{eqnarray}
\left(\frac{\theta}{v}\right)_x&=&\frac{\theta_x}{v}-\frac{\theta v_x}{v^2}\nonumber\\
&\leq& C(v,\theta)(|\theta_x|+|v_x|),\label{2.5}
\end{eqnarray}
\begin{eqnarray}
\left(\frac{\theta}{v}\right)_{xx}&=&
\frac{\theta_{xx}}{v}-\frac{2\theta_xv_x+\theta v_{xx}}{v^2}+\frac{2\theta v_x^2}{v^3}\label{2.6}\\
&\leq& C(v,\theta)\left(\left|\left(v_{xx},\theta_{xx}\right)\right|+\left|\left(v_x,\theta_x\right)\right|^2
\right)\nonumber,
\end{eqnarray}
\begin{eqnarray}
\left(\frac{\theta}{v}\right)_{xxx}&=&
\frac{\theta_{xxx}}{v}-\frac{3\theta_{xx}v_x+3\theta_x v_{xx}+\theta v_{xxx}}{v^2}
+\frac{6\theta_x v_x^2+6\theta v_xv_{xx}}{v^3}-\frac{6\theta v_x^3}{v^4}\label{2.7}\\
&\leq& C(v,\theta)\left(\left|\left(v_x,\theta_x\right)\right|^3+\left|\left(v_{xx},\theta_{xx}\right)\right|
\left|\left(v_{x},\theta_{x}\right)\right|+\left|\left(v_{xxx},\theta_{xxx}\right)\right|\right).\nonumber
\end{eqnarray}
To deduce the energy type estimates on $v(t,x)$, we need some identities on the derivative of $\frac{\mu(\theta)v_x}{v}$ with respect to $x$ up to the second order which are listed below
\begin{eqnarray}
\left(\frac{\mu(\theta)v_x}{v}\right)_x&=&
 \frac{\mu'(\theta)\theta_xv_x+\mu(\theta)v_{xx}}{v}-\frac{\mu(\theta)v_x^2}{v^2}\label{2.8}\\
 &\leq& C(v,\theta)\left(\left|\left(v_x,\theta_x\right)\right|^2+\left|v_{xx}\right|\right),\nonumber
 \end{eqnarray}
 \begin{eqnarray}
\left(\frac{\mu(\theta)v_x}{v}\right)_{xx}-\frac{\mu(\theta)v_{xxx}}{v}&=&
\frac{\mu''(\theta)\theta^2_xv_x+\mu'(\theta)\theta_{xx}v_x
+2\mu'(\theta)\theta_xv_{xx}}{v} \nonumber \\
&&-\frac{2\mu'(\theta)\theta_xv_x^2+3\mu(\theta)v_xv_{xx}}{v^2}+\frac{2\mu(\theta)v_x^3}{v^3}\label{2.9}\\
&\leq& C(v,\theta)\left(\left|\left(v_x,\theta_x\right)\right|^3
+\left|\left(v_x,\theta_x\right)\right|\left|\left(v_{xx},\theta_{xx}\right)\right|
\right)\nonumber.
\end{eqnarray}
Moreover, for derivatives of both $\frac{\mu(\theta)u_x^2}{v}$ and $\frac{\mu(\theta)u_x}{v}$ with respect to $x$ up to the second order or the third order respectively, we have
\begin{eqnarray}\label{2.10}
 \left(\frac{\mu(\theta)u_x^2}{v}\right)_x&=&
 \frac{\mu'(\theta)\theta_xu_x^2+2\mu(\theta)u_xu_{xx}}{v}-\frac{\mu(\theta)u_x^2v_x}{v^2}\\
 &\leq& C(v,\theta)\left(\left|u_xu_{xx}\right|+u_x^2\left(|\theta_x|+|v_x|\right)\right),\nonumber
 \end{eqnarray}
 \begin{eqnarray}\label{2.11}
\left(\frac{\mu(\theta)u_x^2}{v}\right)_{xx}
&=&\frac{\mu''(\theta)\theta_x^2u_x^2+\mu'(\theta)\theta_{xx}u_x^2+4\mu'(\theta)\theta_xu_xu_{xx}
+2\mu(\theta)u^2_{xx}+2\mu(\theta)u_xu_{xxx}}{v} \nonumber \\
&&-\frac{2\mu'(\theta)\theta_xu^2_xv_x+4\mu(\theta)u_xu_{xx}v_x+\mu(\theta)u_x^2v_{xx}}{v^2}
+\frac{2\mu(\theta)u_x^2v_x^2}{v^3}\\
&\leq& C(v,\theta)\left(\left|\left(v_x, u_x, \theta_x\right)\right|^4+u_{xx}^2+
\left|\Big(v_x, u_x, \theta_x\right)\right|^2\left|\left(v_{xx}, u_{xx}, \theta_{xx}\right)\right|\nonumber\\
&&+u_x^2\left|v_{xx}\right|
+\left|u_xu_{xxx}\right|\Big)\nonumber.
\end{eqnarray}
and
\begin{eqnarray}\label{2.12}
\left(\frac{\mu(\theta)u_x}{v}\right)_x&=&\frac{\mu'(\theta)\theta_xu_x+\mu(\theta)u_{xx}}{v}
-\frac{\mu(\theta)u_xv_x}{v^2}\\
&\leq& C(\theta,v)\left(\left|u_{xx}\right|+\left|\left(v_x, u_x,\theta_x\right)\right|^2\right)\nonumber,
\end{eqnarray}
\begin{eqnarray}\label{2.13}
&&\left(\frac{\mu(\theta)u_x}{v}\right)_{xx}-\frac{\mu(\theta)u_{xxx}}{v}\nonumber\\
&=&\frac{\mu''(\theta)\theta_x^2u_x+\mu'(\theta)\theta_{xx}u_x
   +2\mu'(\theta)\theta_xu_{xx}}{v}\nonumber \\
   &&-\frac{2\mu'(\theta)\theta_xv_xu_x+2\mu(\theta)u_{xx}v_x
   +\mu(\theta)u_xv_{xx}}{v^2}+\frac{2\mu(\theta)u_xv^2_x}{v^3}\\
   &\leq& C(\theta,v)\left(\left|\left(v_x,u_x,\theta_x\right)\right|^3+|u_x|\left|v_{xx}\right|
   +\left|\left(v_x,u_x,\theta_x\right)\right|\left|\left(u_{xx},\theta_{xx}\right)\right|
   \right),\nonumber
 \end{eqnarray}
 \begin{eqnarray}\label{2.14}
&&\left(\frac{\mu(\theta)u_x}{v}\right)_{xxx}-\frac{\mu(\theta)u_{xxxx}}{v}\nonumber\\&=&
\frac{\mu'''(v)\theta^3_xv_x+3\mu''(\theta)\theta^2_xu_{xx}+3\mu''(\theta)u_x\theta_x\theta_{xx}
+3\mu'(\theta)\theta_xu_{xxx}}{v}\nonumber \\
&&+\frac{\mu'(\theta)\theta_{xxx}u_{x}+3\mu'(\theta)\theta_{xx}u_{xx}}{v}\nonumber\\
&&-\frac{3\mu''(v)\theta^2_xu_xv_x+3\mu'(\theta)\theta_{xx}u_{x}v_x
+6\mu'(\theta)u_{xx}\theta_xv_x}{v^2}\\
&&-\frac{3\mu'(\theta)\theta_xu_xv_{xx}+3\mu(\theta)u_{xxx}v_{x}+3\mu(\theta)v_{xx}u_{xx}+\mu(\theta)u_xv_{xxx}
}{v^2} \nonumber \\
&&+\frac{6\mu'(v)\theta_xu_xv^2_x+6\mu(v)u_{xx}v_x^2+6\mu(v)v_{xx}u_xv_x}{v^3}
-\frac{6\mu(\theta)u_xv^3_x}{v^4}.\nonumber\\
&\leq& C(v,\theta)\left(\left|\left(v_x,u_x,\theta_x\right)\right|^4
+\left|\left(v_x,u_x,\theta_x\right)\right|^2\left|\left(v_{xx},u_{xx},\theta_{xx}\right)\right|
\right.\nonumber\\
&&+\left|\left(v_x,u_x,\theta_x\right)\right|\left|\left(u_{xxx},\theta_{xxx}\right)\right|+|u_x|\left|v_{xxx}\right|+
\left|u_{xx}\right|\left|\left(v_{xx},\theta_{xx}\right)\right|\Big).\nonumber
 \end{eqnarray}
Finally, for the derivatives of $\frac{\kappa(\theta)\theta_x}{v}$ with respect to $x$ up to the third order, we have
 \begin{eqnarray}
 \left(\frac{\kappa(\theta)\theta_x}{v}\right)_x&=&
 \frac{\kappa'(\theta)\theta^2_x+\kappa(\theta)\theta_{xx}}{v}
 -\frac{\kappa(\theta)\theta_xv_x}{v^2}\label{2.15}\\
 &\leq& C(v,\theta)\left(\theta_x^2+v_x^2+|\theta_{xx}|\right),\nonumber
 \end{eqnarray}
 \begin{eqnarray}
  \left(\frac{\kappa(\theta)\theta_x}{v}\right)_{xx}-\frac{\kappa(\theta)\theta_{xxx}}{v}&=&
 \frac{\kappa''(\theta)\theta^3_x+3\kappa'(\theta)\theta_x\theta_{xx}}{v}
 \nonumber\\
&&-\frac{2\kappa(\theta)\theta_{xx}v_x+2\kappa'(\theta)\theta^2_xv_{x}
 +\kappa(\theta)\theta_xv_{xx}}{v^2}+
 \frac{2\kappa(\theta)\theta_xv^2_x}{v^3},\label{2.16}\\
 &\leq& C(v,\theta)\left(\left|\left(v_x,\theta_x\right)\right|^3
 +\left|\left(v_x,\theta_x\right)\right|\left|\theta_{xx}\right|
 +\left|\theta_x\right|\left|v_{xx}\right|
 \right),\nonumber
 \end{eqnarray}
 \begin{eqnarray}
  &&\left(\frac{\kappa(\theta)\theta_x}{v}\right)_{xxx}-\frac{\kappa(\theta)\theta_{xxxx}}{v}\nonumber\\
  &=&\frac{\kappa'''(\theta)\theta^4_x+6\kappa''(\theta)\theta^2_x\theta_{xx}
  +3\kappa'(\theta)\theta^2_{xx}+4\kappa'(\theta)\theta_x\theta_{xxx}}{v}
 \nonumber\\
 &&-\frac{3\kappa''(\theta)\theta^3_{x}v_x+9\kappa'(\theta)v_x\theta_x\theta_{xx}
 +3\kappa'(\theta)\theta^2_xv_{xx}}{v^2}\label{2.17}\\
 &&-\frac{3\kappa(\theta)\theta_{xxx}v_x+3\kappa(\theta)\theta_{xx}v_{xx}
 +\kappa(\theta)\theta_xv_{xxx}}{v^2}
 \nonumber\\
 &&+\frac{6\kappa(\theta)\theta_{xx}v^2_x+6\kappa'(\theta)\theta^2_xv^2_x
 +6\kappa(\theta)\theta_xv_xv_{xx}}{v^3}-\frac{6\kappa(\theta)\theta_xv^3_x}{v^4}\nonumber\\
 &\leq& C(v,\theta)\left(
 \left|\left(v_x,\theta_x\right)\right|^4+\left|\left(v_x,\theta_x\right)\right|^2
 \left|\left(v_{xx},\theta_{xx}\right)\right|
 +\left|\theta_{xx}\right|\left|\left(v_{xx},\theta_{xx}\right)\right|\right.\nonumber\\
&&+\left|\left(v_x,\theta_x\right)\right|\left|\theta_{xxx}\right|
+\left|\theta_x\right|\left|v_{xxx}\right|\Big).\nonumber
\end{eqnarray}

\section{Energy estimates}
\setcounter{equation}{0}
To prove Theorem 1.1, we first define the following function space for  the solution to the Cauchy problem (\ref{1.5}), (\ref{1.6})
\begin{equation}\label{3.1}
X^k(0,T; M_0,M_1;N_0,N_1)=\left\{\Big(v, u, \theta\Big)(t,x) \left|
\begin{array}{c}
\Big(v-1, u, \theta-1\Big)(t,x)\in C^0\left(0,T;H^k({\bf R})\right)\\[2mm]
\Big(u_x,\theta_x\Big)(t,x)\in L^2\left(0,T; H^k({\bf R})\right)\\[2mm]
M_0\leq v(t,x)\leq M_1,\quad N_0\leq \theta(t,x)\leq N_1.
\end{array}
\right.
\right\}.
\end{equation}
Here $k\geq 1$ is an integer, $T>0$ is a given constant and $M_i,\ N_i$ $(i=0,1)$ are some positive constants.

Under the assumptions given in Theorems 1.1, we can get the following local existence result.

\begin{Lemma} [Local existence] Under the assumptions listed in
Theorem 1.1, there exists a sufficiently small positive constant $t_1$, which depends only
on $\left\|\left(v_0-1,u_0,\frac{\theta_0-1}{\sqrt{\gamma-1}}\right)\right\|_3$, $\underline{V}_0,$ $\overline{V}_0,$  $\underline{\Theta}_0,$ and $\overline{\Theta}_0$, such that the Cauchy problem (\ref{1.5}), (\ref{1.6}) admits a unique smooth solution $(v(t,x), u(t,x), \theta(t,x))$ $
\in X^3\left(0,t_1; \frac 12\underline{V}_0, 2\overline{V}_0; \frac 12\underline{\Theta}_0, 2\overline{\Theta}_0\right)$ and
$(v(t,x), u(t,x), \theta(t,x))$ satisfies
\begin{equation}\label{3.2}
\left\{
\begin{array}{l}
0<\frac 12{\underline{V}_0}\leq v(t,x)\leq 2\overline{V}_0,\\[3mm]
0<\frac 12{\underline{\Theta}_0}\leq \theta(t,x)\leq 2\overline{\Theta}_0
\end{array}
\right.
\end{equation}
for all $(t,x)\in[0,t_1]\times{\bf R}$ and
\begin{equation}\label{3.3}
\max\limits_{t\in[0,t_1]}\left\{\left\|\left(v-1, u, \frac{\theta-1}{\sqrt{\gamma-1}}\right)(t)\right\|_3\right\}\leq 2\left\|\left(v_0-1, u_0,
\frac{\theta_0-1}{\sqrt{\gamma-1}}\right)\right\|_3.
\end{equation}
\end{Lemma}
Lemma 3.1 can be proved by employing the standard iteration argument as in \cite{Antontsev-Kazhikhov-Monakhov, Itaya, Matsumura-Nishida}, the only difference here is that since $\gamma-1$ is sufficiently small in our case, we need to pay particular attention to deal with those terms containing negative powers of $\gamma-1$. Since the modification is straightforward, we thus omit the details for brevity.
\begin{Remark} In Lemma 3.1 the time interval on which the local solution is constructed is claimed to depend on $\left\|\left(v_0-1,u_0,\frac{\theta_0-1}{\sqrt{\gamma-1}}\right)\right\|_3$, an advantage of such a dependence is that we can deduce the estimate (\ref{3.3}) by the smallness of $t_1$. In fact, even if $t_1$ is assumed to depend on $\left\|\left(v_0-1,u_0,\frac{\theta_0-1}{\sqrt{\gamma-1}}\right)\right\|_1$ only, a similar local solvability result of the Cauchy problem (\ref{1.5}), (\ref{1.6}) still holds but in such a case, the local solution $(v(t,x), u(t,x),\theta(t,x))$ constructed in such a way satisfies
\begin{eqnarray}\label{3.4}
\max\limits_{t\in[0,t_1]}\left\{\left\|\left(v-1, u, \frac{\theta-1}{\sqrt{\gamma-1}}\right)(t)\right\|_1\right\}&\leq& 2\left\|\left(v_0-1, u_0,
\frac{\theta_0-1}{\sqrt{\gamma-1}}\right)\right\|_1,\\
\max\limits_{t\in[0,t_1]}\left\{\left\|\left(v-1, u, \frac{\theta-1}{\sqrt{\gamma-1}}\right)(t)\right\|_3\right\}&\leq& C(t_1)\left\|\left(v_0-1, u_0,
\frac{\theta_0-1}{\sqrt{\gamma-1}}\right)\right\|_3.\nonumber
\end{eqnarray}
Here $C(t_1)$ is some positive constant depending only on $t_1$.  Since if we combine the continuation argument with the latter local existence result to extend the local solutions step by step to a global one, the presentation will be rather complex and this is the very reason why we use the local existence result stated in Lemma 3.1.
\end{Remark}

Suppose the local solution $(v(t,x), u(t,x),\theta(t,x))$ constructed in Lemma 3.1 has been extended to the time step $t=T>0$ and satisfies the following {\it a priori} assumptions
\begin{equation}\label{3.5}
\|(\theta-1)(t)\|_3\leq \varepsilon,\quad
0<M^{-1}_1\leq v(t,x)\leq M_1,\quad
\|(v-1,u)(t)\|_3\leq N_1
\end{equation}
for all $x\in{\bf R},\ 0\leq t\leq T$, we now turn to deduce certain energy type estimates on $(v(t,x), u(t,x),$ $\theta(t,x))$ in terms of the initial perturbation. Our main idea here is to use the smallness of both $\varepsilon$ and $\gamma-1$ to control the possible growth of the solution $(v(t,x), u(t,x),\theta(t,x))$ constructed in Lemma 3.1 which is caused by the nonlinearities of the system (\ref{1.5}) under consideration.

In fact, under the assumption that $0<\varepsilon<\min\left\{\overline{\Theta}_0-1,1-\underline{\Theta}_0\right\}$, we have from the {\it a priori} assumption (\ref{3.5}) that
\begin{equation}\label{3.6}
\underline{\Theta}_0\leq\theta(t,x)\leq \overline{\Theta}_0,\quad \|\theta_x(t)\|_{W^{1,\infty}({\bf R})}\leq \varepsilon,\quad \|(u,v)(t)\|_{W^{2,\infty}({\bf R})}\leq N_1+1
\end{equation}
hold for all $x\in{\bf R},\ 0\leq t\leq T$. Without loss of generality, we may assume in the rest of this manuscript that $M_1\geq 1, N_1\geq 1$.

Now we turn to deduce certain energy type estimates on $(v(t,x), u(t,x),\theta(t,x))$. Before doing so, recall that we will use $C$ or $O(1)$ to denote some generic positive constant independent of $\gamma-1$, $M_1$, and $N_1$, but may only depend on the initial data and $C(\cdot,\cdot)$ stands for some positive constant which depends only on the quantities listed in the parenthesis.

The first one is concerned with the basic energy estimate. For this purpose, recall that $R=1$, then it is well-known that
$$
\eta(v,u,\theta)=\phi(v)+\frac{u^2}{2}+\frac{
\phi(\theta)}{\gamma-1}, \quad \phi(x)=x-\ln x-1
$$
is a convex entropy to (\ref{1.5}) which satisfies
$$
\eta(v,u,\theta)_t+\left\{\left(\frac{\theta}{v}-1\right)u
-\frac{\mu(\theta)uu_x}{v}+\frac{(\theta-1)\kappa(\theta)\theta_x}{v\theta}\right\}_x+
\left\{\frac{\mu(\theta)u^2_x}{v\theta}+\frac{\kappa(\theta)\theta_x^2}{v\theta^{2}}\right\}=0.
$$
Integrating the above identity with respect to $t,x$ over $[0,t]\times\bf R$, we have
\begin{Lemma} [Basic energy estimate]
Under the conditions listed in Lemma 3.1, suppose that the local solution $(v(t,x), u(t,x),\theta(t,x))$ constructed in Lemma 3.1 has been extended to the time step $t=T$, then we have for $0\leq t\leq T$ that
\begin{equation}\label{3.7}
\int_{\bf R}\eta(v,u,\theta)dx
+\int^t_0\int_{\bf R}\left(\frac{\mu(\theta)u^2_x}{v\theta}+\frac{\kappa(\theta)\theta_x^2}{v\theta^{2}} \right)dxd\tau=\int_{\bf R}\eta(v_0,u_0,\theta_0)dx.
\end{equation}
\end{Lemma}

Under the {\it a priori} assumption (\ref{3.5}) and (\ref{3.6}), we have from the assumption (\ref{1.4}) and the fact that the upper and lower bounds of $ v_0(x),\theta_0(x)$ do not depend on $\gamma-1$ that
\begin{equation}\label{3.8}
\left\|\left(\sqrt{\phi(v)},u,\frac{\theta-1}{\gamma-1}\right)(t)\right\|^2
+\int_0^t\int_{\bf R} \frac{u_x^2+\theta_x^2}{v}dxd\tau
\leq O(1)\left\|\left(\frac{\theta_0-1}{\sqrt{\gamma-1}},v_0-1,u_0\right)\right\|^2.
\end{equation}

Now we turn to derive the lower and upper bounds on the specific volume $v(t,x)$.
To do so, we need to deduce an estimate on $\left\|\frac{\mu(\theta)v_x}{v}\right\|$.
For this purpose, we have by multiplying (\ref{2.1}) by $\frac{\mu(\theta)v_x}{v}$ that
\begin{eqnarray}
\frac{1}{2}\left(\frac{\mu(\theta)v_x}{v}\right)^2_t+\frac{\mu(\theta)\theta v^2_x}{v^3}
&=&\left(\frac{\mu(\theta)v_xu}{v}\right)_t-\left(\frac{\mu(\theta)v_xu}{v}\right)_x
+\frac{\mu(\theta)u^2_x}{v}\nonumber\\
&&+\frac{\mu(\theta)v_x\theta_x}{v^2}+
 \frac{\mu'(\theta)(u_x\theta_x-v_x\theta_t)(uv-\mu(\theta)v_x)}{v^2}.\label{3.9}
\end{eqnarray}
Integrating the above identity with respect to $t$ and $x$ over $[0,t]\times \bf R$, it follows that
\begin{eqnarray}
&&\frac12\left\|\frac{\mu(\theta)v_x}{v}\right\|^2
+\int_0^t\int_{\bf R}\frac{\mu(\theta)\theta v^2_x}{v^3}dxd\tau\nonumber\\
&=&\frac12\left\|\frac{\mu(\theta_0)v_{0x}}{v_0}\right\|^2
+\underbrace{\int_{\bf R} \left(\frac{\mu(\theta)uv_x}{v}-\frac{\mu(\theta_0)u_0v_{0x}}{v_0}\right)dx}_{I_1}
+\underbrace{\int_0^t\int_{\bf R}\frac{\mu(\theta)u^2_x}{v}dxd\tau}_{I_2}\label{3.10}\\
&&+\underbrace{\int_0^t\int_{\bf R}\frac{\mu(\theta)v_x\theta_x}{v^2}}_{I_3}dxd\tau
+\underbrace{\int_0^t\int_{\bf R}\frac{\mu'(\theta)(u_x\theta_x-v_x\theta_t)(uv-\mu(\theta)v_x)}{v^2}dxd\tau}_{I_4}\nonumber.
\end{eqnarray}
By making use of the Cauchy inequality, the {\it a priori} assumption (\ref{3.5}) and its consequence (\ref{3.6}), and the estimate (\ref{3.8}), it follows that
\begin{eqnarray}
I_1&\leq&\frac14\left\|\frac{\mu(\theta)v_x}{v}\right\|^2
+\|u\|^2+O(1)\left\|(v_{0x},u_0)\right\|^2\nonumber \\
&\leq&\frac14\left\|\frac{\mu(\theta)v_x}{v}\right\|^2
+O(1)\left\|\left(v_0-1,u_0,v_{0x},\frac{\theta_0-1}{\sqrt{\gamma-1}}\right)\right\|^2,\nonumber\\
I_2&\leq& O(1) \left\|\left(v_0-1,u_0,\frac{\theta_0-1}{\sqrt{\gamma-1}}\right)\right\|^2,\label{3.11}\\
I_3&\leq&\frac14\int_0^t\int_{\bf R}\frac{\mu(\theta)\theta v^2_x}{v^3}dxd\tau
+\int_0^t\int_{\bf R}\frac{\mu(\theta)\theta^2_x}{v\theta}dxd\tau\nonumber\\
&\leq&\frac14\int_0^t\int_{\bf R}\frac{\mu(\theta)\theta v^2_x}{v^3}dxd\tau
+O(1)\left\|\left(v_0-1,u_0,\frac{\theta_0-1}{\sqrt{\gamma-1}}\right)\right\|^2\nonumber.
\end{eqnarray}
As to $I_4$, noticing that
\[
u_x\theta_x-v_x\theta_t=u_x\theta_x-(\gamma-1)v_x
\left(\frac{\kappa'(\theta)\theta^2_x+\kappa(\theta)\theta_{xx}+\mu(\theta)u^2_x-u_x\theta}{v}-
\frac{\kappa(\theta)\theta_xv_x}{v^2}\right),
\]
we can get that
\begin{eqnarray}
I_4
&=&\underbrace{\int_0^t\int_{\bf R}\frac{\mu'(\theta)uu_x\theta_x}{v}dxd\tau}_{K_1}
-\underbrace{\int_0^t\int_{\bf R}\frac{\mu(\theta)\mu'(\theta)v_xu_x\theta_x}{v^2}dxd\tau}_{K_2}\nonumber\\
&&-\underbrace{\int_0^t\int_{\bf R}
\frac{(\gamma-1)\mu'(\theta)\left(\kappa'(\theta)u\theta^2_xv_x
+\kappa(\theta)u\theta_{xx}v_x
+\mu(\theta)uv_xu^2_x-u\theta v_xu_{x}\right)}{v^2}dxd\tau}_{K_3}\nonumber\\
&&+\underbrace{\int_0^t\int_{\bf R}\frac{(\gamma-1)\mu'(\theta)v^2_x\left(\mu(\theta)
\kappa'(\theta)\theta^2_{x}+\mu(\theta)\kappa(\theta)\theta_{xx}
+\kappa(\theta)u\theta_{x}+\mu^2(\theta)u^2_{x}-
\mu(\theta)\theta u_{x}\right)}{v^3}dxd\tau}_{K_4}\nonumber\\
&&-\underbrace{\int_0^t\int_{\bf R}
\frac{(\gamma-1)\mu(\theta)\mu'(\theta)\kappa(\theta)\theta_{x}v^3_x}{v^4}dxd\tau}_{K_5}.\label{3.12}
\end{eqnarray}
For $K_i (i=3,4,5)$, we will use the smallness of $\gamma-1$ to control the possible growth of the solutions of the Navier-Stokes equations caused by the nonlinearities of the equations under our consideration. In fact, for $K_4$ and $K_5$, under the {\it a priori} assumption (\ref{3.5}), we can get that
\begin{eqnarray}
K_4+K_5
&=&(\gamma-1)\int_0^t\int_{\bf R}\frac{\mu(\theta)\theta v^2_x}{v^3}\left(
\frac{\mu'(\theta)\kappa'(\theta)\theta^2_x+\mu'(\theta)\kappa(\theta)\theta_{xx}
+\mu'(\theta)\mu(\theta)u^2_x}{\theta}\right.\nonumber\\
&&~~~~~~~~~~\left.+\frac{\mu'(\theta)\kappa(\theta) u\theta_x}{\mu(\theta)\theta}
-\frac{\mu'(\theta)\kappa(\theta) v_x\theta_x}{v\theta}
-\mu'(\theta)u_x\right)dxd\tau\label{3.13}\\
&\leq&O(1)(\gamma-1)\left\|\left(\theta_x^2,\theta_{xx},u\theta_x,v_x^2,\frac{\theta_xv_x}{v},
   u^2_x,u_x\right)\right\|_{L^\infty}\int_0^t\int_{\bf R}\frac{\mu(\theta)\theta v^2_x}{v^3}dxd\tau
\nonumber\\
&\leq& O(1)(\gamma-1)N_1^2M_1\int_0^t\int_{\bf R}\frac{\mu(\theta)\theta v^2_x}{v^3}dxd\tau,\nonumber
\end{eqnarray}
while for $K_3$, we have from the Cauchy inequality, the {\it a priori} assumption (\ref{3.5}), (\ref{3.6}), and the estimate (\ref{3.8}) that
\begin{eqnarray}
|K_3|&\leq&\frac18\int_0^t\int_{\bf R}\frac{\mu(\theta)\theta v^2_x}{v^3}dxd\tau
+O(1)(\gamma-1)^2\int_0^t\int_{\bf R}
\frac{u^2(\theta_x^4+\theta^2_{xx}+u_x^4+u^2_x)}{v}dxd\tau\label{3.14}\\
&\leq&\frac18\int_0^t\int_{\bf R}\frac{\mu(\theta)\theta v^2_x}{v^3}dxd\tau
+O(1)(\gamma-1)^2\left\|\left(u^2\theta_x^2,u^2,u^2u_x^2\right)\right\|_{L^\infty([0,T]\times{\bf R})}\int_0^t\int_{\bf R}\frac{u_x^2+\theta_x^2}{v}dxd\tau\nonumber\\
&&+O(1)(\gamma-1)^2\|u\|^2_{L^\infty([0,T]\times{\bf R})}\int_0^t\int_{\bf R}\frac{\mu(\theta)\theta_{xx}^2}{v}dxd\tau\nonumber\\
&\leq&\frac18\int_0^t\int_{\bf R}\frac{\mu(\theta)\theta v^2_x}{v^3}dxd\tau
+O(1)(\gamma-1)^2\left(\varepsilon^2N_1^2+N_1^4\right)
\left\|\left(v_0-1,u_0,\frac{\theta_0-1}{\sqrt{\gamma-1}}\right)\right\|^2\nonumber\\
&&+O(1)(\gamma-1)^2N_1^2\int_0^t\int_{\bf R}\frac{\mu(\theta)\theta_{xx}^2}{v}dxd\tau\nonumber.
\end{eqnarray}

For $K_1$ and $K_2$, we will use the smallness of $\varepsilon$ to control the possible growth of the solutions mentioned above. To this end, we can deduce from the {\it a priori} assumption (\ref{3.7}), the basic energy estimate (\ref{3.8}), and the Cauchy inequality that
\begin{eqnarray}
K_2
&\leq&\left(\int_0^t\int_{\bf R}\frac{\mu(\theta)\theta v^2_x}{v^3}dxd\tau\right)^{\frac12}
      \left(\int_0^t\int_{\bf R}\frac{\mu(\theta)\left(\mu'(\theta)\right)^2\theta^2_xu^2_x}{v\theta}dxd\tau\right)^{\frac12}\nonumber\\
&\leq & \frac 18\int_0^t\int_{\bf R}\frac{\mu(\theta)\theta v^2_x}{v^3}dxd\tau
+4\left\|\frac{\mu(\theta)(\mu'(\theta))^2\theta_x^2}{\theta}\right\|_{L^\infty([0,T]\times{\bf R})}\int^t_0\int_{\bf R}\frac{u_x^2}{v}dxd\tau\label{3.15}\\
&\leq&\frac 18\int_0^t\int_{\bf R}\frac{\mu(\theta)\theta v^2_x}{v^3}dxd\tau
+O(1)\varepsilon^2\left\|\left(v_0-1,u_0,\frac{\theta_0-1}{\sqrt{\gamma-1}}\right)\right\|^2\nonumber
\end{eqnarray}
and
\begin{eqnarray}
K_1
&\leq&O(1)\int_0^t \|u\|_{L^\infty({\bf R})}\left\|\frac{u_x}{\sqrt v}\right\|
                   \left\|\frac{\theta_x}{\sqrt v}\right\|d\tau\nonumber\\
&\leq&O(1)\int_0^t\|u\|^{\frac12}\|u_x\|^{\frac12}\left\|\frac{u_x}
      {\sqrt v}\right\|\left\|\frac{\theta_x}{\sqrt v}\right\|d\tau\nonumber\\
&\leq&O(1)\int_0^t\int_{\bf R}\frac{u^2_x}{v}dxd\tau+
\int_0^t\|\sqrt v\|_{L^\infty({\bf R})}\left\|\frac 1{\sqrt v}\right\|_{L^\infty({\bf R})}
\|u\|\left\|\frac{u_x}{\sqrt v}\right\|\|\theta_x\|\left\|\frac{\theta_x}{\sqrt v}\right\|d\tau\nonumber\\
&\leq&O(1)\int_0^t\int_{\bf R}\frac{u^2_x}{v}dxd\tau+\varepsilon M_1\|u\|^2\int_0^t\int_{\bf R}
\frac{u^2_x+\theta^2_x}{v}dxd\tau\label{3.16}\\
&\leq& O(1)\left(1+\varepsilon M_1\right)
   \left\|\left(v_0-1,u_0,\frac{\theta_0-1}{\sqrt{\gamma-1}}\right)\right\|^2.\nonumber
\end{eqnarray}

Inserting (\ref{3.13})-(\ref{3.16}) into (\ref{3.12}), we can get that
\begin{eqnarray}
I_4&\leq & \left(\frac 14+O(1)(\gamma-1)M_1N_1^2\right)\int_0^t\int_{\bf R}\frac{\mu(\theta)\theta v^2_x}{v^3}dxd\tau
+O(1)(\gamma-1)^2N_1^2\int_0^t\int_{\bf R}\frac{\mu(\theta)\theta_{xx}^2}{v}dxd\tau\nonumber\\
&&+O(1)\left(1+\varepsilon^2+\varepsilon M_1+(\gamma-1)^2\left(\varepsilon^2N_1^2+N_1^4\right)\right)
\left\|\left(v_0-1,u_0,\frac{\theta_0-1}{\sqrt{\gamma-1}}\right)\right\|^2.\label{3.17}
\end{eqnarray}

Thus, if we plug (\ref{3.11}), (\ref{3.17}) into (\ref{3.10}), it yields that
\begin{eqnarray}
&&\left\|\frac{\mu(\theta)v_x}{v}\right\|^2+
\left(\frac12-C_1(\gamma-1)N^2_1M_1\right)\int_0^t\int_{\bf R}\frac{\mu(\theta)\theta v^2_x}{v^3}dxd\tau\nonumber\\
&\leq& O(1)\|v_{0x}\|^2+O(1)(\gamma-1)^2N^2_1
\int_0^t\int_{\bf R}\frac{\mu(\theta)\theta^2_{xx}}{v}dxd\tau\label{3.18}\\
&&+O(1) \left(1+\varepsilon M_1+(\gamma-1)^2\left(\varepsilon^2 N^2_1+N^4_1\right)\right)
\left\|\left(v_0-1,u_0,\frac{\theta_0-1}{\sqrt{\gamma-1}}\right)\right\|^2.\nonumber
\end{eqnarray}
Here $C_1$ is some positive constant independent of $t$, $x$, and $\gamma-1$.

Having obtained (\ref{3.18}), if we assume that $\gamma-1$ and $\varepsilon$ are small enough such that
$$
\left\{
\begin{array}{l}
\frac14\leq \frac12-C_1(\gamma-1)N^2_1M_1,\\[2mm]
\varepsilon M_1\leq 1,\\[2mm]
(\gamma-1)^2\left(\varepsilon^2 N^2_1+N^4_1\right)\leq 1,
\end{array}
\right.
\leqno(H_1)
$$
then, the above analysis yields the following result
\begin{Lemma}
Under the conditions listed in Lemma 3.2, if we further assume that $\gamma-1$ and $\varepsilon$ are chosen sufficiently small such that (H$_1$) holds true, then we get
\begin{eqnarray}\label{3.19}
\left\|\frac{\mu(\theta)v_x}{v}\right\|^2+
\int_0^t\int_{\bf R}\frac{\mu(\theta)\theta v^2_x}{v^3}dxd\tau
&\leq& O(1)\left\|\left(v_0-1,v_{0x},u_0,\frac{\theta_0-1}{\sqrt{\gamma-1}}\right)\right\|^2\nonumber\\
&&+(\gamma-1)\int_0^t\int_{\bf R}\frac{\mu(\theta)\theta^2_{xx}}{v}dxd\tau.
\end{eqnarray}
\end{Lemma}

From (\ref{3.19}), it is easy to see that to deduce an estimate on $\left\|\frac{\mu(\theta)v_x}{v}\right\|$, we need to deduce an estimate on $\int_0^t\int_{\bf R}\frac{\mu(\theta)\theta^2_{xx}}{v}dxd\tau$ first. To this end, we can get by differentiating
(\ref{1.5})$_3$ with respect to $x$ once and by multiplying the resulting identity by $\theta_x$
that
\begin{eqnarray}
\frac{C_v}{2}\frac{d}{dt}\left(\theta^2_x\right)
+\frac{\kappa(\theta)\theta^2_{xx}}{v}&=&
\left[\left(\frac{\mu(\theta)u^2_x}{v}-u_xp
+\left(\frac{\kappa(\theta)\theta_x}v\right)_x\theta_x\right)\right]_x\label{3.20}\\
&&+\frac{\theta u_x\theta_{xx}}{v}-\frac{\mu(\theta)u^2_x\theta_{xx}}{v}
-\frac{\kappa'(\theta)\theta^2_x\theta_{xx}}{v}
+\frac{\kappa(\theta)\theta_xv_x\theta_{xx}}{v^2}.\nonumber
\end{eqnarray}

Integrating the above equality with respect to $t$ and $x$ over $[0,t]\times\bf R$,
we get
\begin{eqnarray}
&&\frac{1}{2}\left\|\frac{\theta_x}{\sqrt{\gamma-1}}\right\|^2
+\int_0^t\int_{\bf R}\frac{\kappa(\theta)\theta^2_{xx}}{v}dxd\tau\nonumber\\
&=&
\frac{1}{2}\left\|\frac{\theta_{0x}}{\sqrt{\gamma-1}}\right\|^2
+\underbrace{\int_0^t\int_{\bf R}\frac{\theta u_x\theta_{xx}}{v}dxd\tau}_{I_5}
-\underbrace{\int_0^t\int_{\bf R}\frac{\mu(\theta)u^2_x\theta_{xx}}{v}dxd\tau}_{I_6}\nonumber\\
&&-\underbrace{\int_0^t\int_{\bf R}\frac{\kappa'(\theta)\theta^2_x\theta_{xx}}{v}dxd\tau}_{I_7}
+\underbrace{\int_0^t\int_{\bf R}\frac{\kappa(\theta)\theta_xv_x\theta_{xx}}{v^2}dxd\tau}_{I_8}.\label{3.21}
\end{eqnarray}

From the basic energy estimate (\ref{3.8}), the Cauchy inequality, and the {\it a priori} assumption (\ref{3.5}) and (\ref{3.6}), we can get that
\begin{eqnarray}
|I_5|&\leq& \frac15\int_0^t\int_{\bf R}\frac{\kappa(\theta)\theta^2_{xx}}{v}dxd\tau+
O(1)\left\|\left(v_0-1,u_0,\frac{\theta_0-1}{\sqrt{\gamma-1}}\right)\right\|^2,\label{3.22}\\
|I_6|&\leq& \frac15\int_0^t\int_{\bf R}\frac{\kappa(\theta)\theta^2_{xx}}{v}dxd\tau+
O(1)N^2_1\left\|\left(v_0-1,u_0,\frac{\theta_0-1}{\sqrt{\gamma-1}}\right)\right\|^2,\label{3.23}\\
|I_7|&\leq& \frac15\int_0^t\int_{\bf R}\frac{\kappa(\theta)\theta^2_{xx}}{v}dxd\tau+
O(1)\int_0^t\int_{\bf R}\frac{(\kappa'(\theta))^2\theta^4_x}{v\kappa(\theta)}dxd\tau\nonumber\\
 &\leq&\frac15\int_0^t\int_{\bf R}\frac{\kappa(\theta)\theta^2_{xx}}{v}dxd\tau+
O(1)\varepsilon^2\left\|\left(v_0-1,u_0,\frac{\theta_0-1}{\sqrt{\gamma-1}}\right)\right\|^2,\label{3.24}\\
|I_8|&\leq& \frac15\int_0^t\int_{\bf R}\frac{\kappa(\theta)\theta^2_{xx}}{v}dxd\tau
+O(1)\int_0^t\int_{\bf R}\frac{\kappa(\theta)\theta^2_xv^2_x}{v^3}dxd\tau\nonumber\\
 &\leq&\frac15\int_0^t\int_{\bf R}\frac{\kappa(\theta)\theta^2_{xx}}{v}dxd\tau+
 O(1)\left\|\frac{v_x}{v}\right\|^2_{L^\infty([0,T]\times{\bf R})}
 \left\|\left(v_0-1,u_0,\frac{\theta_0-1}{\sqrt{\gamma-1}}\right)\right\|^2\label{3.25}\\
&\leq&\frac15\int_0^t\int_{\bf R}\frac{\kappa(\theta)\theta^2_{xx}}{v}dxd\tau+
 O(1)M^2_1N^2_1\left\|\left(v_0-1,u_0,\frac{\theta_0-1}{\sqrt{\gamma-1}}\right)\right\|^2.
 \nonumber
\end{eqnarray}

Inserting (\ref{3.21})-(\ref{3.25}) into (\ref{3.20}) yields
\begin{eqnarray}\label{3.26}
\left\|\frac{\theta_x}{\sqrt{\gamma-1}}\right\|^2
+\int_0^t\int_{\bf R}\frac{\kappa(\theta)\theta^2_{xx}}{v}dxd\tau
\leq
O(1) \left(N^2_1+M^2_1N^2_1\right)
  \left\|\left(v_0-1,u_0,\frac{\theta_0-1}{\sqrt{\gamma-1}}\right)\right\|_1^2.
\end{eqnarray}

As a direct consequence of the estimates (\ref{3.19}) and (\ref{3.26}), we can deduce that
\begin{Lemma}
Under the same conditions listed in Lemma 3.3, if $\gamma-1$ is further assumed to be sufficiently small such that
$$
(\gamma-1)M^2_1N^2_1\leq1, \leqno(H_2)
$$
then we arrive at
\begin{equation}\label{3.27}
\left\|\frac{\mu(\theta)v_x}{v}\right\|^2+
\int_0^t\int_{\bf R}\frac{\mu(\theta)\theta v^2_x}{v^3}dxd\tau
\leq O(1)\left\|\left(v_0-1,u_0,\frac{\theta_0-1}{\sqrt{\gamma-1}}\right)\right\|^2_1.
\end{equation}
\end{Lemma}

With the estimates (\ref{3.8}) and (\ref{3.27}) in hand, we now apply Y. Kanel's approach, cf. \cite{Kanel}, to deduce a uniform lower bound and a uniform upper bound for $v(t,x)$ . To this end, set
 \begin{equation}\label{3.28}
 \Psi(v)=\int_1^v\frac{\sqrt{\phi(z)}}{z}dz.
 \end{equation}
Note that there exist positive constants $A_1,A_2$ such that
\begin{equation}\label{3.29}
|\Psi(v)|\geq A_1\left(v^{\frac12}+|\ln v|^{\frac32}\right)-A_2.
\end{equation}

Since
\begin{eqnarray}
|\Psi(v)|&=&\left|\int_{-\infty}^x\Psi(v(t,y))_ydy\right|\nonumber\\
&\leq&\int_{\bf R}\left|\frac{\sqrt{\phi(v)}}{v}v_x\right|dx\nonumber\\
       &\leq&\left\|\sqrt{\phi(v)}\right\|\left\|\frac{v_x}{v}\right\|\label{3.30}\\
 &\leq&O(1)\left\|\left(v_0-1,u_0,\frac{\theta_0-1}{\sqrt{\gamma-1}}\right)\right\|^2_1,\nonumber
\end{eqnarray}
we can deduce the following result
\begin{Lemma}
Under the conditions listed in Lemma 3.4, there is positive constant $V_1$, which may depend only on $\underline{\Theta}_0,\ \overline{\Theta}_0,\ \underline{V}_0,\ \overline{V}_0$, and $\left\|\left(v_0-1,u_0,\frac{\theta_0-1}{\sqrt{\gamma-1}}\right)\right\|^2_1$ but independent of $\gamma-1$, such that
\begin{equation}\label{3.31}
V_1^{-1}\leq v(t,x)\leq V_1, \quad\forall (t,x)\in[0,T]\times \bf R.
\end{equation}
\end{Lemma}

To employ the continuation argument to extend the local solutions step by step to a global one, we need to close the {\it a priori} assumption (\ref{3.5}) listed above and for this purpose, we should derive certain higher order energy type estimates.

Firstly, based on the {\it a priori} assumption (\ref{3.6}) and the lower and upper bounds of $v(t,x)$ obtained in Lemma 3.5, we now turn to derive certain energy type estimates on  $(v_x(t,x), u_x(t,x),$ $\theta_x(t,x))$. To this end, we will first get the estimate on $\|u_x\|$.
 \begin{Lemma}
 Under the same conditions of Lemma 3.5, we can get that
 \begin{equation}\label{3.32}
 \|u_x(t)\|^2+\int_0^t\int_{\bf R}\frac{\mu(\theta)u^2_{xx}}{v}dxd\tau
 \leq C(V_1)
 \left\|\left(v_0-1,u_0,\frac{\theta_0-1}{\sqrt{\gamma-1}}\right)\right\|^6_1.
 \end{equation}
 Here and in the rest of this paper, $C(V_1)$ is used to denote some positive constant depending only on $V_1$.
 \end{Lemma}
{\bf Proof:}
Differentiating (\ref{1.5})$_2$ with respect to $x$ once, multiplying the result by $u_x$, and integrating the final identity with respect to $t$ and $x$ over $[0,t]\times \bf R$, we can get that
\begin{eqnarray}\label{3.33}
&&\frac{1}{2}\|u_x\|^2+\int_0^t\int_{\bf R}
\frac{\mu(\theta)u^2_{xx}}{v}dxd\tau\nonumber\\
&=&
\frac{1}{2}\|u_{0x}\|^2
-\underbrace{\int_0^t\int_{\bf R}\left(\frac{\mu(\theta)}{v}\right)_xu_xu_{xx}dxd\tau}_{I_9}
+\underbrace{\int_0^t\int_{\bf R}\left(\frac{\theta}{v}\right)_xu_{xx}dxd\tau}_{I_{10}}.
\end{eqnarray}

The Cauchy inequality together with the estimates (\ref{3.7}) and (\ref{3.27}) yield
\begin{eqnarray}\label{3.34}
I_{10}&\leq&\frac13\int_0^t\int_{\bf R}\frac{\mu(\theta)u^2_{xx}}{v}dxd\tau+
O(1)\int_0^t\int_{\bf R}\left(\frac{\theta^2_x}{v}+\frac{v^2_x}{v^3}\right)dxd\tau\\
&\leq&\frac13\int_0^t\int_{\bf R}\frac{\mu(\theta)u^2_{xx}}{v}dxd\tau
+C(V_1)\left\|\left(v_0-1,u_0,\frac{\theta_0-1}{\sqrt{\gamma-1}}\right)\right\|^2_1,\nonumber\\
I_{9}&\leq&\frac13\int_0^t\int_{\bf R}\frac{\mu(\theta)u^2_{xx}}{v}dxd\tau+
O(1)\int_0^t\int_{\bf R}\left(u^2_xv^2_x+u^2_x\theta^2_x\right)dxd\tau\nonumber\\
&\leq&\frac13\int_0^t\int_{\bf R}\frac{\mu(\theta)u^2_{xx}}{v}dxd\tau
+O(1)\int_0^t\left(\|u_x\|^2_{L^\infty}\|v_x\|^2+\|\theta_x\|^2_{L^\infty}\|u_x\|^2\right)d\tau\label{3.35}\\
&\leq&\frac13\int_0^t\int_{\bf R}\frac{\mu(\theta)u^2_{xx}}{v}dxd\tau
+C(V_1)\left\|\left(v_0-1,u_0,\frac{\theta_0-1}{\sqrt{\gamma-1}}\right)\right\|^2_1
 \int_0^t\|u_x\|\|u_{xx}\|d\tau\nonumber\\
&\leq&\frac12\int_0^t\int_{\bf R}\frac{\mu(\theta)u^2_{xx}}{v}dxd\tau
+C(V_1)\left\|\left(v_0-1,u_0,\frac{\theta_0-1}{\sqrt{\gamma-1}}\right)\right\|^6_1\nonumber.
\end{eqnarray}
Here we have used the fact that $\|\theta_x\|_{L^\infty([0,T]\times{\bf R})}\leq \varepsilon\leq 1$ when dealing with $I_9$.

(\ref{3.34}) and (\ref{3.35}) together with (\ref{3.33}) imply (\ref{3.32}). This completes the proof of Lemma 3.6.

Combine the estimates on $\|u_x\|$ and $\|v_x\|$ obtained in (\ref{3.27}) and (\ref{3.32}) with (\ref{3.21}), we can get the following result
\begin{Lemma}
Under the same conditions listed in Lemma 3.6,
we have
\begin{equation}\label{3.36}
\left\|\frac{\theta_x}{\sqrt{\gamma-1}}\right\|^2
+\int_0^t\int_{\bf R}\frac{\kappa(\theta)\theta^2_{xx}}{v}dxd\tau
\leq
C(V_1)\left\|\left(v_0-1,u_0,\frac{\theta_0-1}{\sqrt{\gamma-1}}\right)\right\|_1^{10}.
\end{equation}
\end{Lemma}
{\bf Proof:} To prove (\ref{3.36}), we only need to deduce better upper bounds on the terms $I_j\ (j=6,7,8)$ on the right hand side of (\ref{3.21}). Since now we have already obtained the uniform lower and upper bounds on $v(t,x)$ and $\theta(t,x)$, we have from the estimate (\ref{3.32}) that
\begin{eqnarray}
I_6&\leq& \frac15\int_0^t\int_{\bf R}\frac{\kappa(\theta)\theta^2_{xx}}{v}dxd\tau+
O(1)\int_0^t\int_{\bf R}u^4_xdxd\tau\nonumber\\
&\leq&\frac15\int_0^t\int_{\bf R}\frac{\kappa(\theta)\theta^2_{xx}}{v}dxd\tau+
O(1)\int_0^t\|u_x\|_{L^\infty_x}^2\left(\int_{\bf R}u^2_xdx\right)d\tau\label{3.37}\\
&\leq&\frac15\int_0^t\int_{\bf R}\frac{\kappa(\theta)\theta^2_{xx}}{v}dxd\tau+
O(1)\int_0^t\|u_x\|^3\|u_{xx}\|d\tau\nonumber\\
&\leq&\frac15\int_0^t\int_{\bf R}\frac{\kappa(\theta)\theta^2_{xx}}{v}dxd\tau+
C(V_1)
  \left\|\left(v_0-1,u_0,\frac{\theta_0-1}{\sqrt{\gamma-1}}\right)\right\|_1^{10}.\nonumber
\end{eqnarray}

Similarly, noticing that the {\it a priori} assumption (\ref{3.5}) implies
$$
\|\theta_x(t)\|_{L^\infty([0,T]\times{\bf R})}\leq\varepsilon\leq 1,
$$
we can deduce that
\begin{eqnarray}
I_7&\leq& \frac15\int_0^t\int_{\bf R}\frac{\kappa(\theta)\theta^2_{xx}}{v}dxd\tau+
C(V_1)\int_0^t\int_{\bf R}\theta^4_xdxd\tau\label{3.38}\\
&\leq& \frac15\int_0^t\int_{\bf R}\frac{\kappa(\theta)\theta^2_{xx}}{v}dxd\tau+C(V_1)
  \left\|\left(v_0-1,u_0,\frac{\theta_0-1}{\sqrt{\gamma-1}}\right)\right\|_1^{2},\nonumber\\
I_8&\leq& \frac15\int_0^t\int_{\bf R}\frac{\kappa(\theta)\theta^2_{xx}}{v}dxd\tau+
C(V_1)\int_0^t\int_{\bf R}\theta^2_xv^2_xdxd\tau \label{3.39}\\
&\leq& \frac15\int_0^t\int_{\bf R}\frac{\kappa(\theta)\theta^2_{xx}}{v}dxd\tau+C(V_1)
  \left\|\left(v_0-1,u_0,\frac{\theta_0-1}{\sqrt{\gamma-1}}\right)\right\|_1^{2}\nonumber.
\end{eqnarray}

Inserting (\ref{3.22}), (\ref{3.37})-(\ref{3.39})
into (\ref{3.26}), we can get (\ref{3.21}) immediately. This completes the proof of Lemma 3.7.

The energy type estimates obtained in Lemma 3.3-Lemma 3.7 imply that under the {\it a priori}
 assumption (\ref{3.5}) and if we assume that $\gamma-1$ and $\varepsilon>0$ are
  chosen sufficiently small such that (H$_1$) and (H$_2$) hold, then there exist 
  two positive constants $V_1\geq 1$, which may depend only on
   $\left\|\left(v_0-1,u_0,\frac{\theta_0-1}{\sqrt{\gamma-1}},v_{0x}\right)\right\|$,
    $\underline{V}_0, \overline{V}_0, \underline{\Theta}_0$, and $\overline{\Theta}_0$ 
    but independent of $T$ and $\gamma-1$, and $C(V_1)$, which depends only on $V_1$
     but independent of $T>0$, $x$, and  $\gamma-1$, such that the following estimates
\begin{eqnarray}\label{3.40}
V^{-1}_1\leq v(t,x)&\leq&  V_1,\quad (t,x)\in[0,T]\times \bf R,\nonumber\\
\left\|\left(v-1,u,\frac{\theta-1}{\sqrt{\gamma-1}}\right)(t)\right\|^2
+\int_0^t\int_{\bf R}\Big(u_x^2&+&\theta_x^2\Big)(\tau,x)dxd\tau\leq C(V_1)\left\|
\left(v_0-1,u_0,\frac{\theta_0-1}{\sqrt{\gamma-1}}\right)\right\|^2,\nonumber\\
\|v_x(t)\|^2+\int_0^t\int_{\bf R}v_x^2(\tau,x)dxd\tau 
&\leq& C(V_1)\left\|\left(v_0-1,u_0,\frac{\theta_0-1}{\sqrt{\gamma-1}}\right)\right\|_1^2\label{3.40},\\
\|u_x(t)\|^2+\int_0^t\int_{\bf R}u_{xx}^2(\tau,x)dxd\tau
 &\leq& C(V_1)\left\|\left(v_0-1,u_0,\frac{\theta_0-1}{\sqrt{\gamma-1}}\right)\right\|_1^6,\nonumber\\
\left\|\frac{\theta_x(\tau)}{\sqrt{\gamma-1}}\right\|^2
+\int_0^t\int_{\bf R}\theta_{xx}^2(\tau,x)dxd\tau
 &\leq& C(V_1)\left\|\left(v_0-1,u_0,\frac{\theta_0-1}{\sqrt{\gamma-1}}\right)\right\|_1^{10}\nonumber
\end{eqnarray}
hold for $0\leq t\leq T$.

Now we turn to derive the second order energy estimates on $(v(t,x), u(t,x), \theta(t,x))$. Firstly for the corresponding estimate on $u_{xx}(t,x)$, we have by
differentiating (\ref{1.5})$_2$ with respect to $x$
twice, multiplying the result by $u_{xx}$, and integrating the final identity with respect to $t$ and $x$ over $[0,t]\times{\bf R}$ that
\begin{eqnarray}\label{3.41}
&&\frac{1}{2}\|u_{xx}(t)\|^2+\int_0^t\int_{\bf R}\frac{\mu(\theta)u^2_{xxx}}{v}dxd\tau\nonumber\\
&=&\frac{1}{2}\|u_{0xx}\|^2+\underbrace{\int_0^t\int_{\bf R}\left(\frac{\theta}{v}\right)_{xx}u_{xxx}dxd\tau}_{I_{11}}\\
&&-\underbrace{\int_0^t\int_{\bf R}\left[\left(\frac{\mu(\theta)u_x}{v}\right)_{xx}
  -\frac{\mu(\theta)u_{xxx}}{v}\right]u_{xxx}dxd\tau}_{I_{12}}\nonumber.
  \end{eqnarray}

(\ref{2.6}) together with the estimate (\ref{3.40}) yield
\begin{eqnarray}
I_{11}&\leq&
\frac14\int_0^t\int_{\bf R}\frac{\mu(\theta)u^2_{xxx}}{v}dxd\tau
+C(V_1)\int_0^t\int_{\bf R}\left(\theta^2_{xx}+\theta^2_xv^2_x+v^2_{xx}+v^4_x\right)dxd\tau\nonumber\\
&\leq&
\frac14\int_0^t\int_{\bf R}\frac{\mu(\theta)u^2_{xxx}}{v}dxd\tau
+C(V_1)\int_0^t\left(\|\theta_{xx}\|^2+
\|\theta_x\|^2\|v_x\|^4+\|v_{xx}\|^2+\|v_{xx}\|\|v_x\|^3\right)d\tau\nonumber\\
&\leq&
\frac14\int_0^t\int_{\bf R}\frac{\mu(\theta)u^2_{xxx}}{v}dxd\tau
+C(V_1)\left\|\left(v_0-1,u_0,\frac{\theta_0-1}{\sqrt{\gamma-1}}\right)\right\|_1^{10}
+C(V_1)\int_0^t\|v_{xx}\|^2d\tau.\label{3.42}
\end{eqnarray}
Similarly, we have from (\ref{2.13}) and (\ref{3.40}) that
\begin{eqnarray}\label{3.43}
I_{12}
&\leq&\frac14\int_0^t\int_{\bf R}\frac{\mu(\theta)u^2_{xxx}}{v}dxd\tau\nonumber\\
&&+C(V_1)\int_0^t\int_{\bf R}\left[\left|\left(u_{xx},\theta_{xx}\right)\right|^2\left|\left(v_x,u_x,\theta_x\right)\right|^2
+v^2_{xx}u^2_x+\left|\left(v_x,u_x,\theta_x\right)\right|^6\right]dxd\tau\nonumber\\
&\leq&\frac14\int_0^t\int_{\bf R}\frac{\mu(\theta)u^2_{xxx}}{v}dxd\tau
+C(V_1)\int_0^t\left\|\left(v_x,u_x,\theta_x\right)\right\|^4_{L^\infty}
\left\|\left(v_x,u_x,\theta_x\right)\right\|^2d\tau\nonumber\\
&&+C(V_1)\int_0^t\int_{\bf R}\left[\left|\left(u_{xx},\theta_{xx}\right)\right|^2\left|\left(v_x,u_x,\theta_x\right)\right|^2
+v^2_{xx}u^2_x\right]dxd\tau\\
&\leq&\frac14\int_0^t\int_{\bf R}\frac{\mu(\theta)u^2_{xxx}}{v}dxd\tau
+C(V_1)\int_0^t\int_{\bf R}\left[\left|\left(u_{xx},\theta_{xx}\right)\right|^2\left|\left(v_x,u_x,\theta_x\right)\right|^2
+v^2_{xx}u^2_x\right]dxd\tau\nonumber\\
&&+C(V_1)\int_0^t\left\|\left(v_x,u_x,\theta_x\right)\right\|^4\left\|\left(v_{xx},u_{xx},
\theta_{xx}\right)\right\|^2d\tau\nonumber\\
&\leq&\frac14\int_0^t\int_{\bf R}\frac{\mu(\theta)u^2_{xxx}}{v}dxd\tau
+C(V_1)\int_0^t\int_{\bf R}\left[\left|\left(u_{xx},\theta_{xx}\right)\right|^2\left|\left(v_x,u_x,\theta_x\right)\right|^2
+v^2_{xx}u^2_x\right]dxd\tau\nonumber\\
&&+C(V_1)\left\|\left(v_0-1,u_0,\frac{\theta_0-1}{\sqrt{\gamma-1}}\right)\right\|_1^{10}
\int_0^t\left\|\left(v_x,u_x,\theta_x\right)\right\|^2\left\|\left(v_{xx},u_{xx},
\theta_{xx}\right)\right\|^2d\tau.\nonumber
\end{eqnarray}

Inserting (\ref{3.42}) and (\ref{3.43}) into (\ref{3.41}), we can deduce that
\begin{eqnarray}\label{3.44}
&&\|u_{xx}\|^2+\int_0^t\|u_{xxx}\|^2d\tau\nonumber\\
&\leq& C(V_1)\|u_{0xx}\|^2+C(V_1)\int_0^t\|v_{xx}\|^2d\tau\nonumber\\
&&+C(V_1)\int_0^t\int_{\bf R}\left[\left|\left(u_{xx},\theta_{xx}\right)\right|^2\left|\left(v_x,u_x,\theta_x\right)\right|^2
+v^2_{xx}u^2_x\right]dxd\tau\\
&&+C(V_1)\left\|\left(v_0-1,u_0,\frac{\theta_0-1}{\sqrt{\gamma-1}}\right)\right\|_1^{10}
\int_0^t\left\|\left(v_x,u_x,\theta_x\right)\right\|^2\left\|\left(v_{xx},u_{xx},
\theta_{xx}\right)\right\|^2d\tau\nonumber\\
&&+C(V_1)\left\|\left(v_0-1,u_0,\frac{\theta_0-1}{\sqrt{\gamma-1}}\right)\right\|_1^{10}.\nonumber
\end{eqnarray}
Here and in the rest of this manuscript, we have assumed without loss of generality that
$$
\left\|\left(v_0-1,u_0,\frac{\theta_0-1}{\sqrt{\gamma-1}}\right)\right\|_1\geq 1.
$$

To deduce an estimate on $\|\theta_{xx}\|$, we have by differentiating (\ref{1.5})$_3$ with respect to $x$ twice, multiplying the result by $\theta_{xx}$, and then integrating the final resulting identity with respect to $t$ and $x$ over $[0,t]\times \bf R$ that
\begin{eqnarray}\label{3.45}
&&\frac12\left\|\frac{\theta_{xx}}{\sqrt{\gamma-1}}\right\|^2+
\int_0^t\int_{\bf R}\frac{\kappa(\theta)\theta^2_{xxx}}{v}dxd\tau\nonumber\\&=&
\frac12\left\|\frac{\theta_{0xx}}{\sqrt{\gamma-1}}\right\|^2+
\underbrace{\int_0^t\int_{\bf R}\left(\frac{\mu(\theta)u^2_x}{v}\right)_x\theta_{xxx}dx\tau}_{I_{13}}+\underbrace{\int_0^t\int_{\bf R}\left(\frac{\theta u_x}{v}\right)_x\theta_{xxx}dxd\tau}_{I_{14}}\nonumber\\
&&\underbrace{-\int_0^t\int_{\bf R}\left(\left(\frac{\kappa(\theta)\theta_x}{v}\right)_{xx}
-\frac{\kappa(\theta)\theta_{xxx}}{v}\right)\theta_{xxx}dx\tau}_{I_{15}}.
\end{eqnarray}

Now we deal with $I_{j},\ (j=13,14,15)$ term by term. For $I_{13}$, we can get from (\ref{2.10}) and (\ref{3.40}) that
\begin{eqnarray}\label{3.46}
I_{13}&\leq& \frac15\int_0^t\int_{\bf R}\frac{\kappa(\theta)}{v}\theta^2_{xxx}dxd\tau
+C(V_1)\int_0^t\int_{\bf R}\left(\theta^2_xu^4_x+u^2_xu^2_{xx}+u^4_xv^2_x\right)dxd\tau\nonumber\\
&\leq&\frac15\int_0^t\int_{\bf R}\frac{\kappa(\theta)}{v}\theta^2_{xxx}dxd\tau
+C(V_1)\int_0^t\int_{\bf R}u^2_xu^2_{xx}dxd\tau\nonumber\\
&&+C(V_1)\int_0^t\|u_x\|^2\|u_{xx}\|^2\left(\|\theta_x\|^2+\|v_x\|^2\right)d\tau\\
&\leq&\frac15\int_0^t\int_{\bf R}\frac{\kappa(\theta)}{v}\theta^2_{xxx}dxd\tau
+C(V_1)\int_0^t\int_{\bf R}u^2_xu^2_{xx}dxd\tau
+C(V_1)\left\|\left(v_0-1,u_0,\frac{\theta_0-1}{\sqrt{\gamma-1}}\right)\right\|_1^{24}.\nonumber
\end{eqnarray}

Moreover, since (\ref{3.40}) implies
\begin{eqnarray*}
&&\int_0^t\int_{\bf R}\left(\theta^4_x\left(\theta^2_x+v^2_x\right)+\theta^2_xv^4_x\right)dxd\tau\\
&\leq&\int_0^t\left(\|\theta_x\|^2\|\theta_{xx}\|^2\left(\|\theta_x\|^2+\|v_x\|^2\right)
+\|v_x\|^2\|\theta_x\|^2\|v_{xx}\|^2\right)d\tau\\
&\leq&C(V_1)\left\|\left(v_0-1,u_0,\frac{\theta_0-1}{\sqrt{\gamma-1}}\right)\right\|_1^{30}+
\int_0^t\|v_x\|^2\|\theta_x\|^2\|v_{xx}\|^2d\tau,
\end{eqnarray*}
we have from (\ref{3.40}), (\ref{2.12}), and (\ref{2.16}) that
\begin{eqnarray}\label{3.47}
I_{15}&\leq&
\frac15\int_0^t\int_{\bf R}\frac{\kappa(\theta)}{v}\theta^2_{xxx}dxd\tau
+C(V_1)\int_0^t\int_{\bf R}\left(\theta^2_x+v^2_x\right)\theta^2_{xx}dxd\tau\nonumber\\
&&+C(V_1)\int_0^t\int_{\bf R}\left(\theta^4_x\left(\theta^2_x+v^2_x\right)+\theta^2_xv^4_x\right)dxd\tau
+C(V_1)\int_0^t\int_{\bf R}\theta^2_xv^2_{xx}dxd\tau\\
&\leq&\frac15\int_0^t\int_{\bf R}\frac{\kappa(\theta)}{v}\theta^2_{xxx}dxd\tau
+C(V_1)\int_0^t\|v_x\|^2\|\theta_x\|^2\|v_{xx}\|^2d\tau\nonumber\\
&&+C(V_1)\left\|\left(v_0-1,u_0,\frac{\theta_0-1}{\sqrt{\gamma-1}}\right)\right\|_1^{30}+
C(V_1)\int_0^t\int_{\bf R}\left(\theta^2_{xx}\left(\theta^2_x+v^2_x\right)+\theta^2_xv^2_{xx}\right)dxd\tau\nonumber
\end{eqnarray}
and
\begin{eqnarray}\label{3.48}
I_{14}&\leq&\frac15\int_0^t\int_{\bf R}\frac{\kappa(\theta)}{v}\theta^2_{xxx}dxd\tau
+C(V_1)\int_0^t\int_{\bf R}\left(\theta^2_xu^2_{x}+u^2_{xx}+u^2_xv^2_x\right)dxd\tau\nonumber\\
&\leq&\frac15\int_0^t\int_{\bf R}\frac{\kappa(\theta)}{v}\theta^2_{xxx}dxd\tau
+C(V_1)\int_0^t\left(\|u_{xx}\|^2+\|\theta_x\|\|u_{x}\|^2\|\theta_{xx}\|+\|u_{xx}\|\|u_x\|\|v_x\|^2\right)d\tau\nonumber\\
&\leq&\frac15\int_0^t\int_{\bf R}\frac{\kappa(\theta)}{v}\theta^2_{xxx}dxd\tau
+C(V_1)\left\|\left(v_0-1,u_0,\frac{\theta_0-1}{\sqrt{\gamma-1}}\right)\right\|_1^{30}.
\end{eqnarray}

Plug (\ref{3.46})-(\ref{3.48}) into (\ref{3.45}), we obtain that
\begin{eqnarray}\label{3.49}
&&\left\|\frac{\theta_{xx}}{\sqrt{\gamma-1}}\right\|^2+
\int_0^t\int_{\bf R}\frac{\kappa(\theta)\theta^2_{xxx}}{v}dxd\tau\nonumber\\
&\leq& C\left\|\frac{\theta_{0xx}}{\sqrt{\gamma-1}}\right\|^2
+C(V_1)\int_0^t\int_{\bf R}\left(\left(u^2_{xx}+\theta^2_{xx}\right)\left(v^2_x+u^2_x+\theta^2_x\right)+\left(u^2_x+\theta_x^2\right)
v^2_{xx}\right)dxd\tau\\
&&+C(V_1)\int_0^t\|v_x\|^2\|\theta_x\|^2\|v_{xx}\|^2d\tau
+C(V_1)\left\|\left(v_0-1,u_0,\frac{\theta_0-1}{\sqrt{\gamma-1}}\right)\right\|_1^{30}.\nonumber
\end{eqnarray}

The estimate on $\|v_{xx}\|$ is complex in some sense. To do so, we first differentiate (\ref{2.1}) with respect to $x$ once, multiply the result by $\left(\frac{\mu(\theta)v_x}{v}\right)_x$, and then integrate the final result with respect to $t$ and $x$ over $[0,t]\times{\bf R}$ to deduce that
\begin{eqnarray}\label{3.50}
\frac12\left\|\left(\frac{\mu(\theta)v_x}{v}\right)_x\right\|^2
&=&
\frac12\left\|\left(\frac{\mu(\theta_0)v_{0x}}{v_0}\right)_x\right\|^2
+\underbrace{\int_0^t\int_{\bf R}u_{tx}\left(\frac{\mu(\theta)v_x}{v}\right)_xdxd\tau}_{I_{16}}\nonumber\\
&&+\underbrace{\int_0^t\int_{\bf R}\left(\frac{\mu(\theta)v_x}{v}\right)_x\left(\frac{\theta} v\right)_{xx}dxd\tau}_{I_{17}}\\
&&+\underbrace{\int_0^t\int_{\bf R}\left(\frac{\mu(\theta)v_x}{v}\right)_x
\left(\frac{\mu'(\theta)}v(\theta_tv_x-\theta_xu_x)\right)_xdxd\tau}_{I_{18}}.\nonumber
\end{eqnarray}

To estimate the terms appear on the right hand side of (\ref{3.50}) term by term, we first deduce from the following identity
\begin{eqnarray}\label{3.51}
I_{16}&=&\int_0^t\int_{\bf R}\left(u_{x}\left(\frac{\mu(\theta)v_x}{v}\right)_x\right)_tdxd\tau
      -\int_0^t\int_{\bf R}u_{x}\left(\frac{\mu(\theta)v_x}{v}\right)_{tx}dxd\tau\nonumber\\
     &=&\int_{\bf R}u_{x}\left(\frac{\mu(\theta)v_x}{v}\right)_xdx-\int_{\bf R}
     u_{0x}\left(\frac{\mu(\theta_0)v_{0x}}{v_0}\right)_xdx
      +\int_0^t\int_{\bf R}u_{xx}\left(\frac{\mu(\theta)v_x}{v}\right)_{t}dxd\tau\\
        &=&\int_{\bf R}u_{x}\left(\frac{\mu(\theta)v_x}{v}\right)_xdx-
     \int_0^t\int_{\bf R}u_{0x}\left(\frac{\mu(\theta_0)v_{0x}}{v_0}\right)_xdx\nonumber\\
     && +\int_0^t\int_{\bf R}\left(\frac{\mu'(\theta)\theta_tv_xu_{xx}}v+
      \frac{\mu(\theta)u^2_{xx}}{v}-\frac{\mu(\theta)v_xu_xu_{xx}}{v^2}\right)dxd\tau\nonumber
\end{eqnarray}
and (\ref{3.40}) that
\begin{eqnarray}\label{3.52}
I_{16}&\leq&\frac14\left\|\left(\frac{\mu(\theta)v_x}{v}\right)_x\right\|^2
+C(V_1)\left\|\left(v_0-1,u_0,\frac{\theta_0-1}{\sqrt{\gamma-1}}\right)\right\|_2^{6}\nonumber\\
&&+C(V_1)\int_0^t\int_{\bf R}\left(|\theta_tv_x|+|v_xu_x|\right)|u_{xx}|dxd\tau.
\end{eqnarray}

(\ref{2.2}) together with (\ref{3.40}) yield
\begin{eqnarray}\label{3.54}
&&\int_0^t\int_{\bf R}|\theta_tv_x||u_{xx}|dxd\tau\nonumber\\
&\leq&C(V_1)(\gamma-1)\int_0^t\int_{\bf R}|v_xu_{xx}|\left(u^2_x+|\theta_{xx}|
+\theta^2_x+|\theta_xv_x|+|u_x|\right)dxd\tau\nonumber\\
&\leq&C(V_1)(\gamma-1)\int_0^t\int_{\bf R}\left(u^2_{xx}+v^2_x\theta^2_{xx}+v^2_x(u^4_x
+\theta^4_x)+v^2_xu^2_x+v_x^4\theta_x^2\right)dxd\tau
\nonumber\\
&\leq&C(V_1)(\gamma-1)\int_0^t\left(\|u_{xx}\|^2+\|v_x\|^2\|\theta_x\|^2\|\theta_{xx}\|^2
+\|v_x\|^2\|u_x\|^2\|u_{xx}\|^2+\|v_x\|^2\right)d\tau\\
&&+C(V_1)(\gamma-1)\int_0^t\|v_x\|^4_{L^\infty}\|\theta_x\|^2d\tau
+C(V_1)(\gamma-1)\int_0^t\int_{\bf R}v^2_x\theta^2_{xx}dxd\tau\nonumber\\
&\leq& C(V_1)\left\|\left(v_0-1,u_0,\frac{\theta_0-1}{\sqrt{\gamma-1}}\right)\right\|_1^{22}
+C(V_1)(\gamma-1)\int_0^t\int_{\bf R}v^2_x\theta^2_{xx}dxd\tau\nonumber\\
&&+C(V_1)(\gamma-1)\int_0^t\|v_{xx}\|^2\|v_x\|^2\|\theta_x\|^2d\tau\nonumber\\
&\leq& C(V_1)\left\|\left(v_0-1,u_0,\frac{\theta_0-1}{\sqrt{\gamma-1}}\right)\right\|_1^{22}
+C(V_1)(\gamma-1)\int_0^t\int_{\bf R}v^2_x\theta^2_{xx}dxd\tau\nonumber\\
&&+C(V_1)(\gamma-1)\left\|\left(v_0-1,u_0,\frac{\theta_0-1}{\sqrt{\gamma-1}}\right)\right\|_1^{12}
\int_0^t\|v_{xx}\|^2d\tau,\nonumber
\end{eqnarray}
and
\begin{eqnarray}\label{3.55}
\int_0^t\int_{\bf R}|u_xv_x||u_{xx}|dxd\tau&\leq&
\int_0^t\int_{\bf R}u^2_{xx}dxd\tau+\int_0^t\int_{\bf R}v^2_xu^2_xdxd\tau\nonumber\\
&\leq&\int_0^t\int_{\bf R}u^2_{xx}dxd\tau+\int_0^t\|u_x\|\|u_{xx}\|\|v_x\|^2dxd\tau\\
&\leq&C(V_1)\left\|\left(v_0-1,u_0,\frac{\theta_0-1}{\sqrt{\gamma-1}}\right)\right\|_1^{6}.\nonumber
\end{eqnarray}
Inserting the estimates of (\ref{3.54}) and (\ref{3.55}) into (\ref{3.52}), it yields that
\begin{eqnarray}\label{3.56}
I_{16}&\leq&
\frac14\left\|\left(\frac{\mu(\theta)v_x}{v}\right)_x\right\|^2
+C(V_1)\left\|\left(v_0-1,u_0,\frac{\theta_0-1}{\sqrt{\gamma-1}}\right)\right\|_1^{22}\\
&&+C(V_1)(\gamma-1)\int_0^t\int_{\bf R}\theta^2_{xx}v^2_xdxd\tau+C(V_1)(\gamma-1)\left\|\left(v_0-1,u_0,\frac{\theta_0-1}{\sqrt{\gamma-1}}\right)\right\|_1^{12}
\int_0^t\|v_{xx}\|^2d\tau.\nonumber
\end{eqnarray}

As to I$_{17}$, we have from (\ref{2.6}), (\ref{2.8}), and (\ref{3.40}) that
\begin{eqnarray}\label{3.57}
I_{17}&=&\int_0^t\int_{\bf R}\left(\frac{\mu'(\theta)\theta_xv_x+\mu(\theta)v_{xx}}{v}-
\frac{\mu(\theta)v^2_x}{v}\right)\left(\frac{\theta_{xx}}{v}-\frac{2\theta_xv_x+\theta v_{xx}}{v^2}
+\frac{2\theta v^2_x}{v^3}\right)dxd\tau\nonumber\\
&\leq&
-\int_0^t\int_{\bf R}\frac{\mu(\theta)\theta v^2_{xx}}{v^3}dxd\tau+
  C(V_1)\int_0^t\int_{\bf R}\frac{\mu(\theta)\theta_{xx} v_{xx}}{v^2}dxd\tau\nonumber\\
&&+C(V_1)\int_0^t\int_{\bf R}\left(|v_{xx}|+|\theta_{xx}|\right)\left(|\theta_xv_x|+v^2_x\right)dxd\tau
+C(V_1)\int_0^t\int_{\bf R}\left(|\theta_xv_x|+v^2_x\right)^2dxd\tau\nonumber\\
&\leq&
-\frac{2}{3}\int_0^t\int_{\bf R}\frac{\mu(\theta)\theta v^2_{xx}}{v^3}dxd\tau+
+C(V_1)\int_0^t\int_{\bf R}\left(\theta^2_{xx}+\theta^2_xv^2_x+v^4_x\right)dxd\tau\\
&\leq&
-\frac12\int_0^t\int_{\bf R}\frac{\mu(\theta)\theta v^2_{xx}}{v^3}dxd\tau+
C(V_1)\int_0^t\left(\|\theta_{xx}\|^2+\|\theta_x\|^2\|v_x\|^4+\|v_x\|^6\right)d\tau\nonumber\\
&\leq&
-\frac12\int_0^t\int_{\bf R}\frac{\mu(\theta)\theta v^2_{xx}}{v^3}dxd\tau+
C(V_1)\left\|\left(v_0-1,u_0,\frac{\theta_0-1}{\sqrt{\gamma-1}}\right)\right\|_1^{10}.\nonumber
\end{eqnarray}

Finally for $I_{18}$, noticing the identities (\ref{2.2}) and (\ref{2.3}) and due to
\begin{eqnarray*}
I_{18}&=&\int_0^t\int_{\bf R}\left(\frac{\mu'(\theta)\theta_xv_x+\mu(\theta)v_{xx}}{v}
-\frac{\mu(\theta)v_x^2}{v^2}\right)
\left(\frac{\mu''(\theta)\theta_x}{v}-\frac{\mu'(\theta)v_x}{v^2}\right)(\theta_tv_x-\theta_xu_x)\\
&&+\frac{\mu'(\theta)}{v}\left(\frac{\mu'(\theta)\theta_xv_x+\mu(\theta)v_{xx}}{v}
-\frac{\mu(\theta)v^2_x}{v^2}\right)
(\theta_{tx}v_x+\theta_tv_{xx}-\theta_{xx}u_x-\theta_xu_{xx})dxd\tau,\nonumber
\end{eqnarray*}
we have
\begin{eqnarray}\label{3.60}
I_{18}&\leq&
C(V_1)\int_0^t\int_{\bf R}\left(|v_x\theta_x|+|v_{x}|^2+|v_{xx}|\right)
\Big\{\left(|\theta_x|+|v_x|\right)\nonumber\\
&&~~~~\times\left[|\theta_xu_x|+(\gamma-1)\left(
|v_xu^2_x|+|v_x\theta_{xx}|+|v_x\theta^2_x|+|v_x|^3\right)\right]\nonumber\\
&&+(\gamma-1)\left(|v_xu_xu_{xx}|+|v_x\theta^3_x|+|v_x\theta_x\theta_{xx}|+
|v_x\theta_{xxx}|+|v_xu_x\theta_x|+|u_x^2\theta_xv_x|\right.\nonumber\\
&&~~~~~~~~~\left.+|v_xu_{xx}|+|v_xu_x|^2+|v_x\theta_x|^2+|v^2_x\theta_{xx}|
   +|v_x\theta_xv_{xx}|+|v^2_xu_x|+|v^3_x\theta_x|\right)\\
 &&+(\gamma-1)\left(|u^2_xv_{xx}|+|v_{xx}\theta_{xx}|+|\theta^2_xv_{xx}|+|v^2_xv_{xx}|
+|u_xv_{xx}|\right)+
 |u_x\theta_{xx}|+|\theta_xu_{xx}| \Big\}dxd\tau.\nonumber\\
&\leq &L_1+L_2+L_3+L_4+L_5.\nonumber
\end{eqnarray}
Here
\begin{eqnarray*}
L_1&=&C(V_1)\int_0^t\int_{\bf R}\left(v^2_x+\theta^2_x+|v_{xx}|\right)
\left(\left|\theta^2_xu_x\right|+\left|v_x\theta_xu_x\right|+\left|\theta_{xx}u_x\right|
+\left|\theta_xu_{xx}\right|\right)dxd\tau,\\
L_2&=&C(V_1)(\gamma-1)\int_0^t\int_{\bf R}\left(v^2_x+\theta^2_x+|v_{xx}|\right)
  \left(\theta^4_x+v^4_x+u^4_x\right)dxd\tau,\\
L_3&=&C(V_1)(\gamma-1)\int_0^t\int_{\bf R}\left(v^2_x+\theta^2_x+|v_{xx}|\right)\left(v^2_x+\theta^2_x+u_x^2\right)
  \left(|\theta_{xx}|+|v_{xx}|+|u_{xx}|\right)dxd\tau,\\
L_4&=&C(V_1)(\gamma-1)\int_0^t\int_{\bf R}\left(v^2_x+\theta^2_x+|v_{xx}|\right)
\left(|\theta_{xx}v_{xx}|+|v_x\theta_{xxx}|\right)dxd\tau,\nonumber\\
L_5&=&C(V_1)(\gamma-1)\int_0^t\int_{\bf R}\left(v^2_x+\theta^2_x+|v_{xx}|\right)
\left(|v_xu_x\theta_x|+|v_xu_{xx}|+|v_x^2u_x|+|u_xv_{xx}|\right)dxd\tau.\nonumber
\end{eqnarray*}

For any $\eta>0$, we have from (\ref{3.40}) that
\begin{eqnarray}\label{3.61}
&&\int_0^t\int_{\bf R}
|v_{xx}|(|u_x|+|\theta_{x}|)(|\theta_{xx}|+|u_{xx}|)dxd\tau\\
&\leq&\eta\int_0^t\int_{\bf R}\frac{\mu(\theta)\theta v^2_{xx}}{v^3}dxd\tau
+ C(\eta)\int_0^t\int_{\bf R}\left(u^2_x+\theta^2_x\right)\left(|\theta^2_{xx}|+|u^2_{xx}|\right)dxd\tau,\nonumber
\end{eqnarray}
\begin{eqnarray}\label{3.62}
&&\int_0^t\int_{\bf R}\left(\theta_x^4+v_x^4+u_x^4\right)dxd\tau\nonumber\\
&\leq&\int_0^t\left(\|\theta_x\|^2_{L^\infty}+\|u_x\|^2_{L^\infty}+\|v_x\|^2_{L^\infty}\right)
\|(\theta_x,u_x,v_x)\|^2d\tau\nonumber\\
&\leq&C(V_1)\left\|\left(v_0-1,u_0,\frac{\theta_0-1}{\sqrt{\gamma-1}}\right)\right\|_1^{10}
\int_0^t\left(\|\theta_x\|\|\theta_{xx}\|+\|u_x\|\|u_{xx}\|+\|v_x\|\|v_{xx}\|\right)d\tau\nonumber\\
&\leq&\eta\int_0^t\int_{\bf R}\frac{\mu(\theta)\theta v^2_{xx}}{v^3}dxd\tau
+C(V_1,\eta)\left\|\left(v_0-1,u_0,\frac{\theta_0-1}{\sqrt{\gamma-1}}\right)\right\|_1^{20}
\int_0^t\|v_x\|^2d\tau\nonumber\\
&&+C(V_1,\eta)\left\|\left(v_0-1,u_0,\frac{\theta_0-1}{\sqrt{\gamma-1}}\right)\right\|_1^{10}
\left(\int_0^t\|\theta_x\|^2d\tau\right)^{\frac12}
\left(\int_0^t\|\theta_{xx}\|^2d\tau\right)^{\frac12}\\
&&+C(V_1,\eta)\left\|\left(v_0-1,u_0,\frac{\theta_0-1}{\sqrt{\gamma-1}}\right)\right\|_1^{10}
\left(\int_0^t\|u_x\|^2d\tau\right)^{\frac12}
\left(\int_0^t\|u_{xx}\|^2d\tau\right)^{\frac12}\nonumber\\
&\leq&\eta\int_0^t\int_{\bf R}\frac{\mu(\theta)\theta v^2_{xx}}{v^3}dxd\tau
+C(V_1,\eta)\left\|\left(v_0-1,u_0,\frac{\theta_0-1}{\sqrt{\gamma-1}}\right)\right\|_1^{22}\nonumber,
\end{eqnarray}
and
\begin{eqnarray}\label{3.63}
\int_0^t\int_{\bf R}\left(\theta_x^4u_x^2+u_x^2v_x^2\theta_x^2\right)dxd\tau
&\leq& \int_0^t\left(\|\theta_x\|^2\|u_x\|^2\|\theta_{xx}\|^2+
  \|v_x\|^2\|u_x\|\|u_{xx}\|\|\theta_x\|\|\theta_{xx}\|\right)d\tau\nonumber\\
&\leq&C(V_1)\left\|\left(v_0-1,u_0,\frac{\theta_0-1}{\sqrt{\gamma-1}}\right)\right\|_1^{26}.
\end{eqnarray}
Consequently we can deduce from (\ref{3.61})-(\ref{3.63}) that
\begin{eqnarray}\label{3.64}
L_1&\leq&\frac{1}{10}\int_0^t\int_{\bf R}\frac{\mu(\theta)\theta v^2_{xx}}{v^3}dxd\tau+
  C(V_1)\left\|\left(v_0-1,u_0,\frac{\theta_0-1}{\sqrt{\gamma-1}}\right)\right\|_1^{26}\nonumber\\
&&+C(V_1)\int_0^t\int_{\bf R}\left(\theta^2_x+v^2_x+u^2_x\right)
  \left(\theta^2_{xx}+u^2_{xx}\right)dxd\tau.
\end{eqnarray}

Repeating the above argument and under the assumption that $\varepsilon$ and $\gamma-1$ are chosen sufficiently small such that (H)$_1$ and (H)$_2$ hold, we can get from (\ref{3.40}) that
\begin{eqnarray}\label{3.65}
L_2&\leq&\frac{1}{20}\int_0^t\int_{\bf R}\frac{\mu(\theta)\theta v^2_{xx}}{v^3}dxd\tau
+C(V_1)(\gamma-1)^2\left(1+N^2_1\right)\int_0^t\int_{\bf R}\left(\theta^6_x+v^6_x+u^6_x\right)dxd\tau\nonumber\\
&\leq&\frac{1}{20}\int_0^t\int_{\bf R}\frac{\mu(\theta)\theta v^2_{xx}}{v^3}dxd\tau
+C(V_1)\int_0^t\left(\|\theta_{xx}\|^2\|\theta_{x}\|^4
  +\|v_{xx}\|^2\|v_{x}\|^4+\|u_{xx}\|^2\|u_{x}\|^4\right)d\tau\nonumber\\
 &\leq&\frac{1}{20}\int_0^t\int_{\bf R}\frac{\mu(\theta)\theta v^2_{xx}}{v^3}dxd\tau
+C(V_1)\left\|\left(v_0-1,u_0,\frac{\theta_0-1}{\sqrt{\gamma-1}}\right)\right\|_1^{30}\\
&&+C(V_1)\left\|\left(v_0-1,u_0,\frac{\theta_0-1}{\sqrt{\gamma-1}}\right)\right\|_1^{2}
\int_0^t\|v_x\|^2\|v_{xx}\|^2d\tau,\nonumber
\end{eqnarray}
\begin{eqnarray}\label{3.66}
L_3&\leq&(\gamma-1)N^2_1\int_0^t\int_{\bf R}|v_{xx}|
(|v_{xx}|+|u_{xx}|+|\theta_{xx}|)dxd\tau
+\frac{(\gamma-1)N^2_1}{20}\int_0^t\int_{\bf R}\frac{\mu(\theta)\theta v^2_{xx}}{v^3}dxd\tau\nonumber\\
&&+C(V_1)(\gamma-1)\int_0^t\int_{\bf R}\left(\theta^2_x+v^2_x+u^2_x\right)^2\left(|u_{xx}|+|\theta_{xx}|\right)dxd\tau\nonumber\\
&&+C(V_1)(\gamma-1)N^2_1\int_0^t\int_{\bf R}\left(\theta^4_x+v^4_x+u^4_x\right)dxd\tau\\
&\leq&C(V_1)(\gamma-1)N^2_1\int_0^t\int_{\bf R}\frac{\mu(\theta)\theta v^2_{xx}}{v^3}dxd\tau
+C(V_1)\left\|\left(v_0-1,u_0,\frac{\theta_0-1}{\sqrt{\gamma-1}}\right)\right\|_1^{22},\nonumber
\end{eqnarray}
\begin{eqnarray}\label{3.67}
L_4&\leq& \frac{\gamma-1}{20}\int_0^t\int_{\bf R}\frac{\mu(\theta)\theta v^2_{xx}}{v^3}dxd\tau
+C(V_1)(\gamma-1)N^2_1\int_0^t\int_{\bf R}\left(|\theta_{xx}v_{xx}|+|\theta_{xxx}v_x|\right)dxd\tau\nonumber\\
&&+C(V_1)
 (\gamma-1)\int_0^t\int_{\bf R}|v_xv_{xx}\theta_{xxx}|dxd\tau\nonumber\\
 &\leq&
 \frac{\gamma-1}{10}\int_0^t\int_{\bf R}\frac{\mu(\theta)\theta v^2_{xx}}{v^3}dxd\tau+C(V_1)(\gamma-1)N^2_1\int_0^t\int_{\bf R}\theta^2_{xxx}dxd\tau\\
 &&+C(V_1)\left(1+(\gamma-1)N^2_1\right)\left\|\left(v_0-1,u_0,\frac{\theta_0-1}{\sqrt{\gamma-1}}\right)\right\|_1^{10},\nonumber
\end{eqnarray}
and
\begin{eqnarray}\label{3.67*}
L_5&\leq& C(V_1)(\gamma-1)\int_0^t\int_{\bf R}\left(v_x^4+\theta_x^4
                 +\theta_x^2v_x^2u_x^2+\theta_x^4u_x+v_x^4u_x^2\right)dxd\tau\nonumber\\
   &&+C(V_1)(\gamma-1)\int_0^t\int_{\bf R}\left(v_{xx}^2+v_x^2u_{xx}^2+u_x^2v_{xx}^2\right)dxd\tau\nonumber\\
   &\leq&\left(C(V_1)(\gamma-1)N_1^2+\eta\right)\int_0^t\int_{\bf R}\frac{\mu(\theta)\theta v^2_{xx}}{v^3}dxd\tau
          +C(V_1)\int_0^t\int_{\bf R}v_x^2u_{xx}^2dxd\tau\nonumber\\
           &&+C(V_1,\eta)\left\|\left(v_0-1,u_0,\frac{\theta_0-1}{\sqrt{\gamma-1}}\right)\right\|_1^{26}.
\end{eqnarray}
Recall that $\eta>0$ is any given sufficiently small positive constant.

Inserting the estimates of (\ref{3.64})-(\ref{3.67*}) into (\ref{3.60}), and if we assume further that
$$
 (\gamma-1)N^4_1\leq 1,\quad  C(V_1)(\gamma-1)N_1^2\leq \frac1{10}
\leqno(H)_3
$$
then, we get
\begin{eqnarray}\label{3.68}
I_{18}&\leq&
C(V_1)\left\|\left(v_0-1,u_0,\frac{\theta_0-1}{\sqrt{\gamma-1}}\right)\right\|_1^{30}
+\frac{1}{5}\int_0^t\int_{\bf R}\frac{\mu(\theta)\theta v^2_{xx}}{v^3}dxd\tau\nonumber\\
&&+C(V_1)\int_0^t\int_{\bf R}\left(\theta^2_x+v^2_x+u^2_x\right)
  \left(\theta^2_{xx}+u^2_{xx}\right)dxd\tau\\
  &&+C(V_1)(\gamma-1)^{\frac 12}\int_0^t\int_{\bf R}\theta^2_{xxx}dxd\tau+
  C(V_1)\left\|\left(v_0-1,u_0,\frac{\theta_0-1}{\sqrt{\gamma-1}}\right)\right\|_1^{2}
\int_0^t\|v_x\|^2\|v_{xx}\|^2d\tau.\nonumber
\end{eqnarray}
Now, inserting (\ref{3.56}), (\ref{3.57}) and (\ref{3.68}) into (\ref{3.50}), we can get that
\begin{eqnarray}\label{3.69}
&&\left\|\left(\frac{\mu(\theta)v_x}{v}\right)_x\right\|^2
+\int_0^t\int_{\bf R}\frac{\mu(\theta)\theta v^2_{xx}}{v^3}dxd\tau\nonumber\\
&\leq&C(V_1)\left\|\left(\frac{\mu(\theta_0)v_{0x}}{v_0}\right)_x\right\|^2
+C(V_1)\int_0^t\int_{\bf R}\left(\theta^2_x+v^2_x+u^2_x\right)
  \left(\theta^2_{xx}+u^2_{xx}\right)dxd\tau\\
&&+C(V_1)\left\|\left(v_0-1,u_0,\frac{\theta_0-1}{\sqrt{\gamma-1}}\right)\right\|_1^{2}
\int_0^t\|v_x\|^2\|v_{xx}\|^2d\tau+C(V_1)(\gamma-1)^{\frac12}\int_0^t\int_{\bf R}\theta^2_{xxx}dxd\tau\nonumber\\
&&+C(V_1)\left\|\left(v_0-1,u_0,\frac{\theta_0-1}{\sqrt{\gamma-1}}\right)\right\|_1^{30}.\nonumber
\end{eqnarray}

Since
\[
\left(\frac{\mu(\theta)v_x}{v}\right)_x=
\frac{\mu'(\theta)\theta_xv_x}{v}+\frac{\mu(\theta)v_{xx}}{v}-\frac{\mu(\theta)v^2_x}{v^2},
\]
we can deduce from (\ref{3.40}) that
\[
\left\|\left(\frac{\mu(\theta)v_x}{v}\right)_x\right\|^2
\geq C(V_1)\|v_{xx}\|^2-C(V_1)\left\|\left(v_0-1,u_0,
\frac{\theta_0-1}{\sqrt{\gamma-1}}\right)\right\|_1^{30}.
\]

Based on the above estimate and (\ref{3.69}), we finally get
\begin{eqnarray}\label{3.70}
&&\left\|v_{xx}\right\|^2
+\int_0^t\|v_{xx}\|^2d\tau\nonumber\\
&\leq&C(V_1)\left\|v_{0xx}\right\|^2
+C(V_1)\int_0^t\int_{\bf R}\left(\theta^2_x+v^2_x+u^2_x\right)
  \left(\theta^2_{xx}+u^2_{xx}\right)dxd\tau\\
&&+C(V_1)\left\|\left(v_0-1,u_0,\frac{\theta_0-1}{\sqrt{\gamma-1}}\right)\right\|_1^{2}
\int_0^t\|v_x\|^2\|v_{xx}\|^2d\tau+C(V_1)(\gamma-1)^{\frac12}\int_0^t\int_{\bf R}\theta^2_{xxx}dxd\tau\nonumber\\
&&+C(V_1)\left\|\left(v_0-1,u_0,\frac{\theta_0-1}{\sqrt{\gamma-1}}\right)\right\|_1^{30}.\nonumber
\end{eqnarray}

A suitable linear combination of (\ref{3.44}), (\ref{3.49}), and (\ref{3.70}) yields the following result
\begin{Lemma}
Under the same condition listed in Lemma 3.7, if we further assume that $\gamma-1$ is 
sufficiently small such that the assumption (H)$_3$ holds, then we have
\begin{eqnarray}\label{3.71}
&&\left\|\left(v_{xx},u_{xx},\frac{\theta_{xx}}{\sqrt{\gamma-1}}\right)\right\|^2
+\int_0^t\left\|\left(v_{xx},u_{xxx},\theta_{xxx}\right)\right\|^2d\tau\nonumber\\
&\leq& C(V_1)\left\|\left(v_0-1,u_0,\frac{\theta_0-1}{\sqrt{\gamma-1}}\right)\right\|_2^{30}
\exp\left(C(V_1)\left\|\left(v_0-1,u_0,\frac{\theta_0-1}{\sqrt{\gamma-1}}\right)\right\|_1^{10}\right).
\end{eqnarray}
\end{Lemma}
{\bf Proof:} In fact, multiplying (\ref{3.70}) by a sufficiently large positive number $\lambda$ and adding the result with (\ref{3.44}) and (\ref{3.49}), we can deduce from the fact that $\gamma-1$ is sufficiently small that
\begin{eqnarray}\label{3.72}
&&\left\|\left(v_{xx},u_{xx},\frac{\theta_{xx}}{\sqrt{\gamma-1}}\right)\right\|^2
+\int_0^t\left\|\left(v_{xx},u_{xxx},\theta_{xxx}\right)\right\|^2d\tau\nonumber\\
&\leq&C(V_1)\int_0^t\int_{\bf R}\left(\theta^2_x+v^2_x+u^2_x\right)
  \left(\theta^2_{xx}+u^2_{xx}\right)dxd\tau
+C(V_1)\int_0^t\int_{\bf R}v^2_{xx}(u^2_x+\theta^2_x)dxd\tau\\ &&+C(V_1)\left\|\left(v_0-1,u_0,\frac{\theta_0-1}{\sqrt{\gamma-1}}\right)\right\|_1^{6}
\int_0^t\left(\|v_x\|^2+\|\theta_x\|^2\right)\left(\|v_{xx}\|^2+\|\theta_{xx}\|^2\right)d\tau\nonumber\\ &&+C(V_1)\left\|\left(v_0-1,u_0,\frac{\theta_0-1}{\sqrt{\gamma-1}}\right)\right\|_1^{30}.\nonumber
\end{eqnarray}

Due to
\begin{eqnarray}
&&\int_0^t\int_{\bf R}v^2_{xx}(u^2_x+\theta^2_x)dxd\tau\nonumber\\
&\leq&
\int_0^t\|v_{xx}\|^2\left(\|u_x\|\|u_{xx}\|+\|\theta_x\|\|\theta_{xx}\|\right)d\tau,\label{3.73}\\
&&\int_0^t\int_{\bf R}\left(\theta^2_x+v^2_x+u^2_x\right)
  \left(\theta^2_{xx}+u^2_{xx}\right)dxd\tau\nonumber\\
  &\leq&
  \int_0^t\left(\|u_{xx}\|^2+\|\theta_{xx}\|^2\right)
    \left(\|u_{x}\|^2+\|u_{xx}\|^2+\|\theta_{x}\|^2+\|\theta_{xx}\|^2\right)d\tau\nonumber\\
&&+ \int_0^t\|v_{x}\|^2
    \left(\|u_{xx}\|\|u_{xxx}\|+\|\theta_{xx}\|\|\theta_{xxx}\|\right)d\tau,\label{3.74}
    \end{eqnarray}
(\ref{3.71}) follows immediately by inserting (\ref{3.73}) and (\ref{3.74})
into (\ref{3.72}) and by employing Gronwall's inequality and the estimate (\ref{3.40}). This completes the proof of Lemma 3.8.

\begin{Remark} Several remarks concerning the second order energy type estimates are given below:
\begin{itemize}
\item[$\bullet$] Since the Navier-Stokes system (\ref{1.5}) is a hyperbolic-parabolic coupled system, the estimates (\ref{3.40}) contain no information on $\int^t_0\|v_{xx}\|^2d\tau$. Fortunately, the term $\int^t_0\int_{\bf R}v_x^2v_{xx}^2dxd\tau$ does not appear on the right hand side of (\ref{3.72}) and consequently our analysis can be continued;
\item[$\bullet$] Since $\theta_t=(\gamma-1)\left(-\frac{R\theta u_x}{v}+\frac{\mu(\theta)u_x^2}{v}+\left(\frac{\kappa(\theta)\theta_x}{v}\right)_x\right)$, to deduce an estimate on $\|v_{xx}\|$, we need to deal with the term $(\gamma-1)\int^t_0\int_{\bf R}\frac{\mu(\theta)\mu'(\theta)\kappa(\theta)}{v^3}\theta_{xx}v^2_{xx}dxd\tau$, cf. (\ref{3.67}) for details. It is easy to see that to bound such a term, we need to deduce an estimate on $\|\theta_{xx}\|_{L^\infty([0,T]\times{\bf R})}$ and as a result, we had to close the energy type estimates in $H^3({\bf R})$.
\end{itemize}
\end{Remark}

Now we deal with the third order energy type estimate on $(v(t,x), u(t,x),\theta(t,x))$. For this purpose, to simplify the presentation, we denote the $H^3({\bf R})-$norm of the initial perturbation by
$$
N_0=\left\|\left(\frac{\theta_0-1}{\sqrt{\gamma-1}},v_0-1,u_0\right)\right\|_3,
$$
then, from (\ref{3.40}) and (\ref{3.71}), we know that there exists a nonnegative smooth function $C(N_0)$ satisfying $C(0)=0$ such that
\begin{equation}\label{3.75}
\left\|\left(\frac{\theta-1}{\sqrt{\gamma-1}},v-1,u\right)(t)\right\|^2_2+
\int_0^t\left(\|v_x(\tau)\|^2_1+\left\|\left(u_x,\theta_x\right)(\tau)\right\|_2^2\right)d\tau
\leq C(N_0).
\end{equation}
(\ref{3.75}) together with Sobolev's imbedding inequality imply
\begin{equation}\label{3.76}
\left\|\left(\frac{\theta-1}{\sqrt{\gamma-1}},v-1,u \right)(t)\right\|_{W^{1,\infty}({\bf R})}
\leq C(N_0),\quad 0\leq t\leq T.
\end{equation}

We now turn to deduce the desired third order energy type estimates on $(v(t,x), u(t,x),\theta(t,x))$. To this end, we first consider the estimate on $\theta_{xxx}$ and obtain by differentiating (\ref{1.5})$_3$ with respect to $x$
three times, multiplying the result by $\theta_{xxx}$, and integrating the final result with respect to $t$ and $x$ over $[0,t]\times{\bf R}$ that
\begin{eqnarray}\label{3.77}
&&\frac12\left\|\frac{{\theta}_{xxx}}{\sqrt{\gamma-1}}\right\|^2
+\int_0^t\int_{\bf R}\frac{\kappa(\theta)\theta^2_{xxxx}}{v}dxd\tau\nonumber\\
&=&\frac12\left\|\frac{{\theta}_{0xxx}}{\sqrt{\gamma-1}}\right\|^2
\underbrace{-\int_0^t\int_{\bf R}\theta_{xxxx}\left(\frac{\mu(\theta)u^2_x}{v}\right)_{xx}dxd\tau}_{I_{19}}\\
&&\underbrace{-\int_0^t\int_{\bf R}\theta_{xxxx}\left(\left(\frac{\kappa(\theta)\theta_x}{v}\right)_{xxx}
  -\frac{\kappa(\theta)\theta_{xxxx}}{v}\right)dxd\tau}_{I_{20}}\nonumber\\
  &&\underbrace{+\int_0^t\int_{\bf R}\theta_{xxxx}\left(\frac{\theta u_x}{v}\right)_{xx}dxd\tau}_{I_{21}}.\nonumber
\end{eqnarray}
To estimate $I_j\ (j=19, 20, 21)$ term by term, we have from (\ref{2.11}), (\ref{3.75}), and (\ref{3.76}) that
\begin{eqnarray}\label{3.78}
I_{19}
&\leq&\frac{1}{10}\int_0^t\int_{\bf R}\frac{\kappa(\theta)\theta^2_{xxxx}}{v}dxd\tau
+C(V_1)\int_0^t\int_{\bf R}\left(\theta^4_xu^4_x+\theta^2_{xx}u^4_x+\theta^2_xu^2_xu^2_{xx}+u^4_{xx}\right.\nonumber\\
&&~~~~~~~~~~~~~~~~~\left.+u^2_{xxx}u^2_x+\theta^2_xv^2_xu^4_x+v^2_{xx}u^4_x+v^4_xu^4_x\right)dxd\tau\\
 &\leq& \frac{1}{10}\int_0^t\int_{\bf R}\frac{\kappa(\theta)\theta^2_{xxxx}}{v}+C(N_0).\nonumber
\end{eqnarray}
Here we have use the fact that
\[
\int_0^t\int_{\bf R}u^4_{xx}dxd\tau\leq\int_0^t\|u_{xxx}\|\|u_{xx}\|^3d\tau
\leq \int_0^t\|u_{xxx}\|^2+\int_0^t\|u_{xx}\|^6d\tau\leq C(N_0).
\]

Similarly, (\ref{2.17}) and (\ref{2.13}) together with (\ref{3.75}) and
(\ref{3.76}) imply
\begin{equation}\label{3.79}
I_{20}
\leq \frac{1}{10}\int_0^t\int_{\bf R}\frac{\kappa(\theta)\theta^2_{xxxx}}{v}dxd\tau
+C(N_0)+C(N_0)\int_0^t\int_{\bf R}\left(\theta^2_{xx}v^2_{xx}+v^2_{xxx}\right)dxd\tau
\end{equation}
and
\begin{equation}\label{3.80}
I_{21}
\leq \frac{1}{10}\int_0^t\int_{\bf R}\frac{\kappa(\theta)\theta^2_{xxxx}}{v}dxd\tau
+C(N_0).
\end{equation}

Inserting (\ref{3.78})-(\ref{3.80}) into (\ref{3.77}) yields
\begin{eqnarray}\label{3.81}
&&\left\|\frac{{\theta}_{xxx}}{\sqrt{\gamma-1}}\right\|^2
+\int_0^t\|\theta_{xxxx}(\tau)\|^2d\tau\nonumber\\
&\leq& C(N_0)
+C(N_0)\int_0^t\int_{\bf R}\left(\theta^2_{xx}v^2_{xx}+v^2_{xxx}\right)dxd\tau\\
&\leq& \left\|\frac{{\theta}_{0xxx}}{\sqrt{\gamma-1}}\right\|^2
+C(N_0)+C(N_0)\int_0^t\|\theta_{xx}\|\|\theta_{xxx}\|\|v_{xx}\|^2d\tau
+C(N_0)\int_0^t\int_{\bf R}v^2_{xxx}dxd\tau\nonumber\\
&\leq& \left\|\frac{{\theta}_{0xxx}}{\sqrt{\gamma-1}}\right\|^2
+C(N_0)+C(N_0)\int_0^t\left(\|\theta_{xx}\|^2+\|v_{xx}\|^2\right)\|\theta_{xxx}\|^2d\tau
+C(N_0)\int_0^t\int_{\bf R}v^2_{xxx}dxd\tau,\nonumber
\end{eqnarray}
and we have from the Gronwall inequality and (\ref{3.75}) and
(\ref{3.76}) that
\begin{equation}\label{3.82}
\left\|\frac{{\theta}_{xxx}}{\sqrt{\gamma-1}}\right\|^2
+\int_0^t\|\theta_{xxxx}(\tau)\|^2d\tau
\leq C(N_0)+C(N_0)\int_0^t\int_{\bf R}v^2_{xxx}dxd\tau.
\end{equation}

To deduce an estimate on $u_{xxx}$, we have by differentiating (\ref{1.5})$_2$ with respect to $x$ three times, multiplying the result by $u_{xxx}$, and then integrating the final result with respect to $t$ and $x$ over $[0,t]\times {\bf R} $ that
\begin{eqnarray}\label{3.83}
&&\frac12\|u_{xxx}\|^2+\int_0^t\int_{\bf R}\frac{\mu(\theta)u^2_{xxxx}}{v}dxd\tau\nonumber\\
&=&\frac12\|u_{0xxx}\|^2
\underbrace{-\int_0^t\int_{\bf R}\left(\left(\frac{\mu(\theta)u_{x}}{v}\right)_{xxx}-\frac{\mu(\theta)u_{xxxx}}{v}\right)
u_{xxxx}dxd\tau}_{I_{22}}\\
&&+
\underbrace{\int_0^t\int_{\bf R}\left(\frac{\theta}{v}\right)_{xxx}u_{xxxx}dxd\tau}_{I_{23}}.\nonumber
\end{eqnarray}

 (\ref{2.7}) and (\ref{2.14}) together with the estimates (\ref{3.75}) and (\ref{3.76}) imply
\begin{eqnarray}\label{3.84}
I_{22}&\leq& \frac15\int_0^t\|u_{xxxx}\|^2d\tau+
 C(N_0)+C(N_0)\int_0^t\int_{\bf R}\left(u^2_{xx}\theta^2_{xx}+u^2_{xx}v^2_{xx}+v^2_{xxx}\right)dxd\tau\nonumber\\
 &\leq&\frac15\int_0^t\|u_{xxxx}\|^2d\tau+
 C(N_0)+C(N_0)\int_0^t\int_{\bf R}v^2_{xxx}dxd\tau
 \end{eqnarray}
and
\begin{equation}\label{3.85}
I_{23}\leq\frac15\int_0^t\|u_{xxxx}\|^2d\tau
+C(N_0)\int_0^t\|v_{xxx}\|^2d\tau+C(N_0).
\end{equation}

Putting (\ref{3.84}),(\ref{3.85}), and (\ref{3.83}) together, we can obtain
\begin{equation}\label{3.86}
\|u_{xxx}\|^2+\int_0^t\|u_{xxxx}\|^2dxd\tau
\leq C(N_0)
+C(N_0)\int_0^t\int_{\bf R}v^2_{xxx}dxd\tau.
\end{equation}

Finally to get an estimate on $v_{xxx}$, we have from (\ref{2.1}) that
\begin{equation}\label{3.87}
\left(\frac{\mu(\theta)v_x}{v}\right)_{txx}
=u_{txx}+\left(\frac{\theta}{v}\right)_{xxx}
+\left(\frac{\mu'(\theta)}{v}\left(\theta_tv_x-u_x\theta_x\right)\right)_{xx}.
\end{equation}
Multiply (\ref{3.87}) by $\left(\frac{\mu(\theta)v_x}{v}\right)_{xx}$
and integrate the resulting identity with respect to $t$ and $x$
over$[0,t]\times \bf R$, we have
\begin{eqnarray}\label{3.88}
\frac{1}{2}\left\|\left(\frac{\mu(\theta)v_x}{v}\right)_{xx}\right\|^2
&=&\underbrace{\int_0^t\int_{\bf R}\left(\frac{\mu(\theta)v_x}{v}\right)_{xx}u_{txx}dxd\tau}_{I_{24}}+
\underbrace{\int_0^t\int_{\bf R}\left(\frac{\mu(\theta)v_x}{v}\right)_{xx}\left(\frac{\theta}{v}\right)_{xxx}dxd\tau}_{I_{25}}\nonumber\\
&&+\underbrace{\int_0^t\int_{\bf R}\left(\frac{\mu(\theta)v_x}{v}\right)_{xx}
\left(\frac{\mu'(\theta)}{v}\left(\theta_tv_x-u_x\theta_x\right)\right)_{xx}dxd\tau}_{I_{26}}.
\end{eqnarray}

To deal with the terms appeared on the right hand side of (\ref{3.88}) term by term, we first have
\begin{eqnarray}\label{3.89}
I_{24}&=&\int_0^t\int_{\bf R}\left[\left(u_{xx}\left(\frac{\mu(\theta)v_x}{v}\right)_{xx}\right)_t
      -u_{xx}\left(\frac{\mu(\theta)v_x}{v}\right)_{txx}\right]dxd\tau\nonumber\\
&=&\int_{\bf R}u_{xx}\left(\frac{\mu(\theta)v_x}{v}\right)_{xx}dx-\int_{\bf R}
      u_{0xx}\left(\frac{\mu(\theta_0)v_{0x}}{v_0}\right)_{xx}dx
      +\int_0^t\int_{\bf R}u_{xxx}\left(\frac{\mu(\theta)v_x}{v}\right)_{tx}dxd\tau\nonumber\\
  &=&\int_{\bf R}u_{xx}\left(\frac{\mu(\theta)v_x}{v}\right)_{xx}dx-\int_{\bf R}
     u_{0xx}\left(\frac{\mu(\theta_0)v_{0x}}{v_0}\right)_{xx}dx\nonumber\\
      &&+\int_0^t\int_{\bf R}\left(\frac{\mu'(\theta)\theta_tv_x}v+
      \frac{\mu(\theta)u_{xx}}{v}-\frac{\mu(\theta)v_xu_x}{v^2}\right)_xu_{xxx}dxd\tau\nonumber\\
   &\leq& C(N_0)+\frac15 \left\|\left(\frac{\mu(\theta)v_x}{v}\right)_{xx}\right\|^2
     +C(N_0)\underbrace{\int_0^t\int_{\bf R}|\theta_{tx}v_xu_{xxx}|dxd\tau}_{J_1}\\
    &&+C(N_0)\underbrace{\int_0^t\int_{\bf R}|u_{xxx}\theta_t|\left(|\theta_xv_x|+|v^2_x|+|v_{xx}|\right)dxd\tau}_{J_2}\nonumber\\
    && +C(N_0)\underbrace{\int_0^t\int_{\bf R}|u_{xxx}|\left(|u_{xxx}|+|v_xu_{xx}|+|\theta_xu_{xx}|+|v^2_xu_{x}|
      +|v_{xx}u_{x}|+|v_xu_{x}\theta_x|
    \right)dxd\tau}_{J_3}\nonumber\\
   &\leq&C(N_0)+\frac15 \left\|\left(\frac{\mu(\theta)v_x}{v}\right)_{xx}\right\|^2.\nonumber
\end{eqnarray}
Here we have used the following estimates
\begin{eqnarray*}
J_1&\leq& C(N_0)\int_0^t\int_{\bf R}\left(v^2_xu^2_{xxx}+\theta^2_{tx}\right)dxd\tau
\leq C(N_0),\\
J_2&\leq& C(N_0)\int_0^t\int_{\bf R}\left(u^2_{xxx}\left(\theta^2_xv^2_x+v^4_x\right)+\theta^2_{t}+\theta^2_{t}v^2_{xx}\right)dxd\tau\nonumber\\
&\leq& C(N_0)+C(N_0)(\gamma-1)^2\int_0^t\int_{\bf R}\theta^2_{xx}v^2_{xx}dxd\tau\\
&\leq&C (N_0)+C (N_0)(\gamma-1)^2\int_0^t\|\theta_{xx}\|\|\theta_{xxx}\|\|v_{xx}\|^2d\tau\nonumber\\
&\leq&C (N_0)+C (N_0)(\gamma-1)^2\int_0^t\|\theta_{xxx}\|^2d\tau\nonumber\\
&\leq&C (N_0),\nonumber\\
J_3&\leq& \int_0^t\|u_{xxx}\|^2d\tau+
   C(N_0)\int_0^t\int_{\bf R}\left(u^2_{xx}v^2_{x}+\theta^2_{x}u^2_{xx}+u^2_{x}v^4_{x}
    +u^2_{x}v^2_{xx}+\theta^2_{x}v^2_{x}u^2_x\right)dxd\tau\nonumber\\
    &\leq& C(N_0)
\end{eqnarray*}
which follow from (\ref{2.2}), (\ref{2.3}), (\ref{3.75}), and (\ref{3.76})

As to the term I$_{25}$, we get from (\ref{2.7}), (\ref{2.9}), (\ref{3.75}) and (\ref{3.76}) that
\begin{eqnarray}\label{3.94}
I_{25}
&\leq& -\int_0^t\int_{\bf R}\frac{\mu(\theta)\theta v^2_{xxx}}{v^3}dxd\tau\nonumber\\
&&+C(N_0)\int_0^t\int_{\bf R}\Bigg\{ \left(|\theta_{xxx}|+|v_{x}\theta_{xx}|+|v_{xx}\theta_{x}|+|v_x|^3
+|v_{x}v_{xx}|+|\theta_xv_x^2|\right)\nonumber\\
&&~~~~~~~~~~~~\times\left(|\theta^2_{x}v_{x}|+|v_{x}\theta_{xx}|
+|\theta_{x}v_{xx}|+|v_xv_{xx}|+|v_x|^3\right)\nonumber\\
&&+|v_{xxx}|\left(|\theta_{xxx}|+|v_{x}\theta_{xx}|+|v_{xx}\theta_{x}|+|v_x|^3
+|v_{x}v_{xx}|+|\theta_xv_x^2|+|\theta^2_xv_x|+|\theta_xv_{xx}|\right)\Bigg\}dxd\tau\nonumber\\
&\leq&-\frac12\int_0^t\int_{\bf R}\frac{\mu(\theta)\theta v^2_{xxx}}{v^3}dxd\tau
 +C(N_0).
\end{eqnarray}

To treat the term $I_{26}$ is much more complex than the other terms on the right hand side of (\ref{3.88}), although this process is similar to the proof of I$_{18}$, we shall give the proof in detail for reader's convenience. In fact, notice that
\begin{eqnarray}
&&\left(\frac{\mu'(\theta)}{v}\left(\theta_tv_x-u_x\theta_x\right)\right)_{x}\label{3.95}\\
&\leq &C(N_0)\left(|\theta_t\theta_xv_x|+|\theta_tv_{xx}|+|\theta_tv^2_x|+
    |\theta_{tx}v_x|+|\theta^2_{x}u_x|+|u_{xx}\theta_x|+|u_x\theta_{xx}|+|\theta_xv_xu_x|
\right),\nonumber\\
&&\left(\frac{\mu'(\theta)}{v}\left(\theta_tv_x-u_x\theta_x\right)\right)_{xx}\nonumber\\
&\leq& C(N_0)\left(|\theta_t\theta^2_xv_x|+|\theta_t\theta_{xx}v_x|+|\theta_t\theta_xv_{xx}|
+|\theta_t\theta_xv^2_x|+|\theta_tv_{xxx}|+|\theta_tv_{xx}v_x|+|\theta_tv^3_x|\right.\label{3.96}\\
&&+|\theta_{tx}\theta_xv_x|+|\theta_{tx}v_{xx}|+|\theta_{tx}v^2_x|+|\theta_{txx}v_{x}|
+|\theta^3_xu_x|+|\theta_xu_x\theta_{xx}|+|\theta^2_xu_{xx}|+|\theta^2_xu_xv_{x}|\nonumber\\
&&+|\theta_xu_{xxx}|+|\theta_xu_{xx}v_{x}|+|u_x\theta_{xxx}|+|\theta_xu_xv_{xx}|
\left. +|u_{xx}\theta_{xx}|+|v_xu_x\theta_{xx}|+|\theta_xu_xv^2_{x}|\right),\nonumber
\end{eqnarray}
we have from (\ref{3.75}),(\ref{3.76}) and noticing that $\gamma-1$ can be chosen sufficiently small that
\begin{eqnarray}\label{3.97}
I_{26}&\leq& C(N_0)+\frac1{10}\int_0^t\int_{\bf R}\frac{\mu(\theta)\theta v^2_{xxx}}{v^3}dxd\tau\nonumber\\
&&+C(N_0)(\gamma-1)^2\int_0^t\int_{\bf R}\left(\theta^2_{xxxx}+\theta^2_{xx}v^2_{xxx}+\theta^2_{xxx}v^2_{xx}\right)dxd\tau\nonumber\\
&\leq& C(N_0)+\frac1{10}\int_0^t\int_{\bf R}\frac{\mu(\theta)\theta v^2_{xxx}}{v^3}dxd\tau\nonumber\\
&&+C(N_0)(\gamma-1)^2\int_0^t\int_{\bf R}\theta^2_{xxxx}dxd\tau
+C(N_0)\int_0^t\|v_{xxx}\|^2\|\theta_{xx}\|\|\theta_{xxx}\|d\tau.
\end{eqnarray}

Now, we insert (\ref{3.89}), (\ref{3.94}), and (\ref{3.97}) into (\ref{3.88}) and
let $\gamma-1$  small enough, we can obtain from the Gronwall inequality that
\begin{eqnarray}\label{3.98}
&&\left\|\left(\frac{\mu(v)v_x}{v}\right)_{xx}\right\|^2
+\int_0^t\int_{\bf R}\frac{\mu(\theta)\theta v^2_{xxx}}{v^3}dxd\tau\\
&\leq &C(N_0)
+C(N_0)(\gamma-1)^2\int_0^t\|\theta_{xxxx}\|^2d\tau d\tau.\nonumber
\end{eqnarray}

Due to
\[
\left\|\left(\frac{\mu(v)v_x}{v}\right)_{xx}\right\|\geq C(N_0)\|v_{xxx}\|^2-C(N_0),
\]
we have from (\ref{3.98}) that
\begin{equation}\label{3.99}
\|v_{xxx}\|^2+\int_0^t\|v_{xxx}\|^2d\tau
\leq C(N_0)+C(N_0)(\gamma-1)^2\int_0^t\|\theta_{xxxx}\|^2d\tau.
\end{equation}

For sufficiently large positive constant $\lambda$, we have by performing (\ref{3.99})$\times\lambda+$(\ref{3.82})$+$(\ref{3.86}) and noticing that $\gamma-1$ can be chosen sufficiently small that
\begin{Lemma}
Under the same conditions listed in Lemma 3.8, if we let $\gamma-1$ small enough, then we have
\begin{equation}\label{3.100}
\left\|\left(v_{xxx},u_{xxx},\frac{\theta_{xxx}}{\sqrt{\gamma-1}}\right)(t)\right\|^2
+\int^t_0\left(\left\|v_{xxx}(\tau)\right\|^2+\left\|\left(u_{xxxx},\theta_{xxxx}\right)(\tau)\right\|^2\right)d\tau\leq C(N_0).
\end{equation}
\end{Lemma}

As a direct consequence of Lemmas 3.3-3.9, we have
\begin{Corollary} [Energy type {\it a priori} estimates] Let $(v(t,x), u(t,x), \theta(t,x))$ be the local solution constructed in Lemma 3.1 which has been extended to the time step $t=T\geq t_1$ and  assume that $(v(t,x), u(t,x), \theta(t,x))$ satisfies the {\it a priori} assumption (\ref{3.5}), then if $\varepsilon>0$ and $\gamma-1>0$ are chosen sufficiently small such that
$$
\left\{
\begin{array}{l}
(\gamma-1)C_2(N_0)<\frac 12,\\[2mm]
\frac 14\leq \frac 12-C_1(\gamma-1)N_1^2M_1,\\[2mm]
(\gamma-1)\left(M_1^4+N_1^4\right)\leq 1,\\[2mm]
C(V_1)(\gamma-1)N_1^2\leq \frac1{10},\\[2mm]
(\gamma-1)^2\left(\varepsilon^2N_1^2+N_1^4\right)\leq 1,\\[2mm]
\varepsilon M_1\leq 1,
\end{array}
\right.
\leqno{(H)}
$$
there exists a nonnegative function $C_3(N_0)$ satisfying $C_3(0)=0$ such that
\begin{equation}\label{3.101}
\left\|\left(v-1,u,\frac{\theta-1}{\sqrt{\gamma-1}}\right)(t)\right\|_3^2
+\int^t_0\left(\left\|v_{x}(\tau)\right\|_2^2+\left\|\left(u_{x},\theta_{x}\right)(\tau)\right\|_3^2\right)d\tau\leq C_3(N_0)
\end{equation}
holds for any $0\leq t\leq T$ and there exists a positive constant which depends only on $N_0$ such that
\begin{equation}\label{3.102}
V_1^{-1}\leq v(t,x)\leq V_1,\quad  \forall(t,x)\in[0,T]\times {\bf R}.
\end{equation}
Here $C_1$ is defined in (\ref{3.18}) and $C_i(N_0)\ (i=2,3)$ are some positive constants depending only on $N_0$.

Moreover if we assume further that $\gamma-1>0$ is sufficiently small such that
\begin{equation}\label{3.103}
(\gamma-1)C_3(N_0)\leq \min\left\{\left(1-\underline{\Theta}_0\right)^2, \left(\overline{\Theta}_0-1\right)^2\right\},
\end{equation}
then we have
\begin{equation}\label{3.104}
\underline{\Theta}_0\leq \theta(t,x)\leq \overline{\Theta}_0,\quad \forall(t,x)\in[0,T]\times {\bf R}.
\end{equation}
\end{Corollary}

\section{The proof of our main result}
\setcounter{equation}{0}

This section is devoted to proving our main result which is based on the continuation argument.
Before doing so, noticing that
$\theta=\frac{A}{R}v^{1-\gamma}\exp\left(\frac{\gamma-1}{R}s\right)$, $\overline s=\frac{R}{\gamma-1}\ln{\frac{R}{A}}$,
and recalling that we have assume that $A=1,\ R=1$ which imply that $\overline{s}=0$, we have
\begin{eqnarray*}
\theta-1&=&v^{1-\gamma}\exp\left((\gamma-1)s\right)-1\\
        &=&v^{1-\gamma}\exp\left((\gamma-1)s\right)-
         \exp\left((\gamma-1)\overline s\right)\\
          &=&\left(v^{1-\gamma}-1\right)\exp\left((\gamma-1)s\right)
          +\exp\left((\gamma-1)s\right)-\exp\left((\gamma-1)
                    \overline{s}\right).
\end{eqnarray*}
Consequently
\begin{eqnarray}\label{4.1}
\|\theta_0-1\|&\leq&O(1)(\gamma-1)\exp\left((\gamma-1)\|s_0\|_{L_x^\infty}\right)
     \left[\left\|v_0^{-\gamma}\right\|_{L_x^\infty}\|v_0-1\|+\|s_0(x)-\overline s\|\right],\nonumber\\
\|\theta_{0x}\|&\leq&O(1)(\gamma-1)\exp\left((\gamma-1)\|s_{0}\|_{L_x^\infty}\right)
                 \left[\left(\inf\limits_x v_0(x)\right)^{-\gamma}\|v_{0x}\|
                 +\left(\inf_x v_0(x)\right)^{1-\gamma}\|s_{0x}\|\right],\nonumber\\
\|\theta_{0xx}\|&\leq&O(1)(\gamma-1)\exp\left((\gamma-1)\|s_{0x}\|_{L_x^\infty}\right)
                 \Big[\left(\inf\limits_x v_0(x)\right)^{-\gamma}\|(v_{0xx},s_{0xx})\|\nonumber\\
           &&+\left(\inf\limits_x v_0(x)\right)^{-\gamma-1}\|(v_{0x}^2,s_{0x}^2)\|\Big],\nonumber\\
 \|\theta_{0xxx}\|&\leq&O(1)(\gamma-1)\exp\left((\gamma-1)\|s_{0x}\|_{L_x^\infty}\right)
  \left[\left(\inf\limits_x v_0(x)\right)^{-\gamma-2}\|(v_{0x}^3,s_{0x}^3)\|\right.\nonumber \\
     &&+\left.\left(\inf\limits_x v_0(x)\right)^{-\gamma-1}\|(v_{0x},s_{0x})(v_{0xx},s_{0xx})\|
      +\left(\inf\limits_x v_0(x)\right)^{-\gamma}\|(v_{0xxx},s_{0xxx})\|\right].
\end{eqnarray}
Thus, if we assume that $\|v_0(x)\|_{L_x^\infty}, \inf\limits_x v_0(x), \frac{\gamma-1}{A}\|s_0(x)\|_{L_x^\infty}$ are independent of $\gamma-1$, we have from (\ref{4.1}) and the assumptions listed in Theorem 1.1 that
\begin{equation}\label{4.2}
\left\|\frac{\theta_0-1}{\sqrt{\gamma-1}}\right\|_3\leq C_4\sqrt{\gamma-1}\left\|\left(v_0-1,u_0,s_0-\overline{s}\right)\right\|_3.
\end{equation}
Here $C_4$ is some positive constant independent of $\gamma-1$.

Now we turn to prove Theorem 1.1. First under the conditions listed in Theorem 1.1, we have from the local existence result stated in Lemma 3.1 that there exists a sufficiently small positive constant $t_1$, which depends only on $\underline{V}_0, \overline{V}_0,  \underline{\Theta}_0, \overline{\Theta}_0$ and $\|(v_0-1,u_0,\theta_0-1)\|_3$, such that the Cauchy problem (\ref{1.5}), (\ref{1.6}) admits a unique smooth solution $(v(t,x), u(t,x), \theta(t,x)
\in X^3\left(0,t_1; \frac 12\underline{V}_0, 2\overline{V}_0; \frac 12\underline{\Theta}_0,
2\overline{\Theta}_0\right)$ which satisfies
\begin{equation}\label{4.3}
\frac 12\underline{V}_0\leq v(t,x)\leq 2\overline{V}_0,\quad \frac 12\underline{\Theta}_0\leq \theta(t,x)\leq 2\overline{\Theta}_0
\end{equation}
and
\begin{equation}\label{4.4}
\left\|\left(v-1, u, \frac{\theta-1}{\sqrt{\gamma-1}}\right)(t)\right\|_3\leq 2\left\|\left(v_0-1, u_0,
\frac{\theta_0-1}{\sqrt{\gamma-1}}\right)\right\|_3:=2N_0
\end{equation}
for all $0\leq t\leq t_1,\ x\in{\bf R}$.

(\ref{4.2}) together with (\ref{4.4}) imply
\begin{equation}\label{4.5}
\|(v-1,u)(t)\|_3\leq 2N_0,\quad \|\theta(t)-1\|_3\leq 2\sqrt{\gamma-1}N_0,\quad 0\leq t\leq t_1.
\end{equation}
Thus if we set $T=t_1$, $M_1=\min\left\{2\overline{V}_0,2\underline{V}_0^{-1}\right\}$, $\varepsilon=2\sqrt{\gamma-1}N_0$, $N_1=2N_0$, since $N_0$, $\overline{V}_0$, $\underline{V}_0$, $\overline{\Theta}_0$, and $\underline{\Theta}_0$ are assumed to be independent of $\gamma-1$,one can easily deduce that there exists a positive constant $\gamma_1>1$ such that if $1<\gamma\leq \gamma_1$ we have that
\begin{equation}\label{4.6}
\left\{
\begin{array}{l}
\left(\gamma-1\right)C_2(N_0)<\frac 12,\\[2mm]
\frac 14\leq \frac 12-C_1\left(\gamma-1\right)N_1^2M_1,\\[2mm]
(\gamma-1)\left(M_1^4+N_1^4\right)\leq 1,\\[2mm]
C(V_1)(\gamma-1)N_1^2\leq \frac1{10},\\[2mm]
(\gamma-1)^2\left(\varepsilon^2N_1^2+N_1^4\right)\leq 1,\\[2mm]
2\sqrt{\gamma-1}N_0 M_1\leq 1
\end{array}
\right.
\end{equation}
and
\begin{equation}\label{4.7}
2\sqrt{\gamma-1}N_0\leq \min\left\{\left(1-\underline{\Theta}_0\right)^2, \left(\overline{\Theta}_0-1\right)^2\right\}
\end{equation}
hold. The above analysis tells us that all the conditions listed in Corollary 3.1 hold with $T=t_1$ and consequently we have from Corollary 3.1 that the local solution $(v(t,x), u(t,x), \theta(t,x))$ constructed above satisfies
\begin{eqnarray}\label{4.8}
&&V_1^{-1}\leq v(t,x)\leq V_1,\quad \underline{\Theta}_0\leq \theta(t,x)\leq \overline{\Theta}_0,\\
&&\left\|\left(v-1,u,\frac{\theta-1}{\sqrt{\gamma-1}}\right)(t)\right\|_3^2
+\int^t_0\left(\left\|v_{x}(\tau)\right\|_2^2+\left\|\left(u_{x},\theta_{x}\right)(\tau)\right\|_3^2\right)d\tau\leq C_3(N_0)\nonumber
\end{eqnarray}
for all $0\leq t\leq t_1,\ x\in{\bf R}.$ Here $V_1$ and $C_3(N_0)$ are some positive constants defined in Corollary 3.1.

Now, we take  $(v(t_1,x), u(t_1,x), \theta(t_1,x))$ as initial data, we have from the estimates (\ref{4.8}) and Lemma 3.1 that the local solution $(v(t,x), u(t,x), \theta(t,x))$ constructed above can be extended to the time step $t=t_1+t_2$ for some suitably small positive constant $t_2$ depending only on $N_0$, $V_1$, $\underline{\Theta}_0$, and $\overline{\Theta}_0$ and satisfies
\begin{equation}\label{4.9}
\frac {1}{2V_1}\leq v(t,x)\leq 2V_1,\quad \frac 12\underline{\Theta}_0\leq \theta(t,x)\leq 2\overline{\Theta}_0
\end{equation}
and
\begin{equation}\label{4.10}
\left\|\left(v-1, u, \frac{\theta-1}{\sqrt{\gamma-1}}\right)(t)\right\|_3\leq 2\left\|\left(v-1, u,
\frac{\theta-1}{\sqrt{\gamma-1}}\right)(t_1)\right\|_3\leq 2\sqrt{C_3(N_0)}
\end{equation}
for all $0\leq t\leq t_1+t_2,\ x\in{\bf R}$. If we set $T=t_1+t_2$, $\varepsilon=2\sqrt{(\gamma-1)C_3(N_0)}$, $M_1=2V_1$, $N_1=2\sqrt{C_3(N_0)}$, since $V_1$ and $N_0$ are independent of $\gamma-1$, it is easy to see that we can find a positive constant $\gamma_2>0$ such that for all $1<\gamma\leq \gamma_2$
\begin{equation}\label{4.11}
\left\{
\begin{array}{l}
\left(\gamma-1\right)C_2(N_0)<\frac 12,\\[2mm]
\frac 14\leq \frac 12-C_1\left(\gamma-1\right)N_1^2M_1,\\[2mm]
(\gamma-1)\left(M_1^4+N_1^4\right)\leq 1,\\[2mm]
C(V_1)(\gamma-1)N_1^2\leq \frac1{10},\\[2mm]
(\gamma-1)^2\left(\varepsilon^2N_1^2+N_1^4\right)\leq 1,\\[2mm]
2\sqrt{\gamma-1}N_0 M_1\leq 1
\end{array}
\right.
\end{equation}
and
\begin{equation}\label{4.12}
2\sqrt{(\gamma-1)C_3(N_0)}\leq \min\left\{\left(1-\underline{\Theta}_0\right)^2, \left(\overline{\Theta}_0-1\right)^2\right\}
\end{equation}
hold. The above analysis tells us that all the conditions listed in Corollary 3.1 hold with $T=t_1+t_2$ and consequently we have from Corollary 3.1 that the solution $(v(t,x), u(t,x), \theta(t,x))$ defined on the time interval $[0,t_1+t_2]$ satisfies (\ref{4.8}) for all $0\leq t\leq t_1+t_2$ with the same positive constants $V_1$ and $C_3(N_0)$.

Now, we take  $(v(t_1+t_2,x), u(t_1+t_2,x), \theta(t_1+t_2,x)$ as initial data, noticing that the constants $V_1$ and $C_3(N_0)$ in (\ref{4.8}) are independent of the time variable $t$, we can then extend  $(v(t,x),$ $ u(t,x), \theta(t,x)$ to the time step $t=t_1+2t_2$ by exploiting Lemma 3.1 again. Repeating the above procedure, if we take $\gamma_0=\min\{\gamma_1,\gamma_2\}$, we can thus extend the solution $(v(t,x), u(t,x), \theta(t,x))$ step by step to a global one provided that
$$
1<\gamma\leq \gamma_0
$$
and as a by-product of the above analysis, we can also deduce that $(v(t,x), u(t,x), \theta(t,x))$ satisfies
\begin{equation}\label{4.13}
\left\|\left(v-1,u,\frac{\theta-1}{\sqrt{\gamma-1}}\right)(t)\right\|_3^2+
\int_0^t\left(\left\|v_x(\tau)\right\|_2^2
+\left\|\left(u_x,\theta_x\right)(\tau)\right\|_3^2\right)d\tau\leq C_3(N_0).
\end{equation}
From which the time asymptotic behavior (\ref{1.10}) follows easily. This completes the proof of Theorem 1.1.

\bigbreak

\begin{center}
{\bf Acknowledgment}
\end{center}
Hongxia Liu was supported by a grant from the National Natural Science Foundation of China under contract 11271160. Tong Yang was supported by Joint Research Fund by National Natural Science Foundation of China and Research Grants Council of Hong Kong, N-CityU102/12.
 Huijiang Zhao was supported by a grant from the National Natural Science Foundation of China under contract 10925103. Qingyang Zou was supported by ``the Fundamental Research Funds for the Central Universities''. This work was also supported by a grant from the National Natural Science Foundation of China under contract 11261160485.

\end{document}